\documentclass[12pt]{article}
\usepackage{amsmath}
\usepackage{amsfonts}
\usepackage{amssymb}
\usepackage{latexsym}
\usepackage{graphicx}
\usepackage{enumerate}
\usepackage[shortlabels]{enumitem}
\usepackage{mathtools}
\usepackage{hyperref}

\usepackage{amsmath, amssymb,latexsym}
\usepackage{color}
\usepackage{xcolor}
\usepackage{fullpage}

\def\C{\mathbb{C}}
\def\F{\mathbb{F}}
\def\la{\lambda}

\def\Z{\mathbb{Z}}
\def\T{\mathsf{T}}

\DeclareMathOperator{\diag}{diag}

\DeclareMathOperator{\rank}{rank}

\newtheorem{theorem}{Theorem}[section]
\newtheorem{proposition}[theorem]{Proposition}

\newtheorem{lemma}[theorem]{Lemma}
\newtheorem{definition}[theorem]{Definition}
\newtheorem{corollary}[theorem]{Corollary}
\newtheorem{remark}[theorem]{{\sc Remark}}
\newtheorem{example}[theorem]{Example}

\title{Root vectors of polynomial and rational matrices: theory and computation}
\author{
Vanni Noferini\thanks{Aalto University, Department of Mathematics and Systems Analysis, P.O. Box 11100, FI-00076, Aalto, Finland. Supported by an Academy of Finland grant (Suomen Akatemian p\"{a}\"{a}t\"{o}s 331240).}
\and 
Paul Van Dooren\thanks{Universit\'e catholique de Louvain, Department of Mathematical Engineering, Av. Lemaitre 4, B-1348 Louvain-la-Neuve, Belgium. Supported by an Aalto Science Institute Visitor Programme.}
}
\begin{document}
\maketitle
\begin{abstract} 
The notion of root polynomials of a polynomial matrix $P(\la)$ was thoroughly studied in [F. Dopico and V. Noferini, Root polynomials and their role in the theory of matrix polynomials, Linear Algebra Appl. 584:37--78, 2020]. In this paper, we extend such a systematic approach to
general rational matrices $R(\lambda)$, possibly singular and possibly with coalescent pole/zero pairs. We discuss the related theory for rational matrices with coefficients in an arbitrary field. As a byproduct, we obtain sensible definitions of eigenvalues and eigenvectors of a rational matrix $R(\la)$, without any need to assume that $R(\la)$ has full column rank or that the eigenvalue is not also a pole. Then, we specialize to the complex field and provide a practical algorithm to compute them, based on the construction of a minimal state space realization of the rational matrix $R(\lambda)$ and then using the staircase algorithm on the linearized 
pencil to compute the null space as well as the root polynomials in a given point $\la_0$. If $\la_0$ is also a pole, then it is necessary to apply a preprocessing step that removes the pole while making it possible to recover the root vectors of the original matrix: in this case, we study both the relevant theory (over a general field) and an algorithmic implementation (over the complex field), still based on minimal state space realizations.
\end{abstract}



\section{Introduction}

Let $R(\la)\in \C(\la)^{m\times n}$ be a rational matrix with coefficients in the field of complex numbers $\C$. In linear system theory and in control theory, an important problem is to find the zeros and poles of $R(\la)$. In this setting, $R(\la)$ contains the information of the inputs and outputs of a system of differential or difference equations \cite{kailath,Rosenbrock}. It is shown in \cite{MFK} that the natural frequencies of the dynamical system are captured by the poles, and the frequencies that are blocked by the system are the zeros. When $R(\la)$ is singular, that it is, either it is not square or it does not have full rank over the field of rational functions, then the matrix has also left and/or right null spaces \cite{For75}. These vectors spaces and their so-called minimal bases provide further information on the initial conditions and on the degrees of freedom associated with the response of the system \cite{For75,kailath,MFK,Rosenbrock}. Additional motivation to study matrix-valued rational functions stems from the nonlinear eigenvalue problem \cite{GT}. In this context, one has a locally meromorphic matrix-valued function $M(\la)$, that could either be exactly rational \cite{GT,MV} or possibly approximated with a rational function $R(\la)$ \cite{GT,GBMM}. One method to approximate (some of the) zeros and eigenvectors of $M(\la)$ is to first compute zeros and eigenvectors of $R(\la)$.

A point $\la_0\in \C$ may be a simple or multiple zero, a simple or multiple pole, and even a coalescent pole and zero of the rational matrix $R(\la)$. The local (at $\la_0$) Smith-McMillan form of the rational matrix is then relevant. In linear systems theory, it refines the description of how the dynamical system responds at the frequency $\la_0$ via the \emph{partial multiplicities} of $\la_0$. In nonlinear eigenvalue problems, the partial multiplicities are also important. They generalize concepts familiar for the classical eigenvalue problem, like the size of the Jordan blocks of a matrix, or (equivalently) the Segr\'{e} and Weyr characteristics of a matrix. In particular, information about the partial multiplicities is highly relevant for the perturbation theory of an eigenvalue $\la_0$.

The spectral description of a rational matrix does not, however, end with the partial multiplicities. In linear systems theory, the response also depends on particular {\it directions}, which are input vectors that {\it excite} the system's pole or zero of a multiplicity given by each structural index.
These vectors arise naturally when one wants to describe the solution set of particular matrix equations involving rational matrices \cite{kailath,kar94} and are also used in tangential interpolation problems of high order \cite{GVV04}. Such vectors can be viewed as a generalization of an eigenvector of a first order system of differential equation modeled by the eigenvalue problem $\la x- Ax=0$.  Similarly, in nonlinear eigenvalue problems one may be interested in appropriate generalizations of the eigenvectors of a matrix. If the nonlinear matrix is regular, and the eigenvalue is not also a pole, then eigenvectors can indeed be always defined as usual. In the case of polynomial eigenvalue problems, eigenvectors were also defined for singular matrix polynomials \cite[Definition 2.18]{DopN21}.

More generally, in this paper we show how to extend the theory of {\it root polynomials} \cite{DopN21, GLR82} to the case of rational matrices over any field. For a regular polynomial (or rational) matrix, root polynomials can be viewed as generating functions for the set of Jordan chains at a point. In the singular case, they generalize such a concept, and as mentioned above they for example can be used to define eigenvectors as elements of certain quotient space, see \cite[Subsection 2.3]{DopN21}. We will show that root polynomials of rational matrices encode all the information on structural indices \cite{DopN21}, and that they make it possible to properly define eigenvectors for rational matrices, even when they are singular or when they have coalescent poles and zeros (or even when they present both issues).

The paper is organized as follows. Sections \ref{Sec:background} and \ref{Sec:rootvectors} are theoretical, and there we work over an algebraically closed field $\F$ (not necessarily $\C$); this is not a restriction because the case of a non-closed field can be dealt with by embedding. In Section \ref{Sec:background} we recall the basic definitions and background material for the rest of the paper, including the concept of root polynomials at a point $\la_0$ for a polynomial matrix. In  Section \ref{Sec:rootvectors} we recall the concept of local discrete valuations and use this to define root vectors of a general rational matrix with a coalescent pole and zero at $\la_0$. We develop a general theory of root vectors at a finite point showing, amongst other things, that they can be assumed without loss of generality to be root polynomials, that they can be used to define eigenvectors even when the rational matrix has a pole at the point of interest (or even when it is singular), and that \emph{maximal sets} of root vectors always encode the information about partial multiplicities. We also mention how to define root vectors at infinity. The rest of the paper is concerned on the computation of root polynomials of polynomial matrices and root vectors of rational matrices. For this we first recall what is a minimal realization of a general rational matrix $R(\la)$ in Section \ref{Sec:systemmatrix}. Then, in Section \ref{Sec:rootsofSandR}, we show how to recover maximal sets of root polynomials of $R(\la)$ from those of a minimal system matrix, first assuming that $\la_0$ is not a pole; the case of $\la_0$ being also a pole is later dealt with by constructing a minimal realization of a second rational matrix $Y(\la)$ that has the property of having no poles at $\la_0$, but whose maximal sets of root vectors at $\la_0$ are simply related to those of $R(\la)$ and allow us to recover them. We then specialize to the case where the system matrix is a minimal realization, and hence a pencil. In the case of $\F\subseteq\C$, algorithms to compute root polynomials of pencils exist \cite{NofV}, and hence overall this yields an algorithm for computing maximal sets of root polynomials of any rational matrix $R(\la)$. Such an algorithm is relatively simple when $\la_0$ is not also a pole, but also works (albeit with additional complications) in the case of coalescent poles/zeros. Finally, we make some concluding remarks in Section \ref{Sec:conclusion}.

\section{Background and definitions} \label{Sec:background}
Let $\F$ be an algebraically closed field. The theory that we will develop can be generalized to any field $\mathbb{K}$ by embedding $\mathbb{K}^{m \times n}$ in $\F^{m \times n}$ where $\F$ is the algebraic closure of $\mathbb{K}$. The assumption is made for notational convenience as it guarantees that any finite eigenvalue is an element of the base field. Given a point $\lambda_0 \in \F$, the local Smith-McMillan form of the $m \times n$ rational matrix
$R(\lambda)\in \F(\lambda)^{m\times n}$ at $\lambda_0$ is
:
\begin{equation} \label{MCM}
M(\lambda)\cdot R(\lambda)\cdot N(\lambda) :=
\left[ \begin{array}{ccc|c} (\lambda-\lambda_0)^{\sigma_r} & & 0 & \\
 & \ddots & & \\ 0 & &  (\lambda-\lambda_0)^{\sigma_1}& \\ \hline
 & & & 0_{m-r,n-r} \end{array}\right].
\end{equation}
In \eqref{MCM}, $M(\lambda)$ and $N(\lambda)$ are rational matrices that do not have poles at $\la_0$ and are invertible at
$\lambda_0$, $r$ is the rank of $R(\lambda)$ over $\F(\la)$ (sometimes called the {\em normal rank} of $R(\lambda)$) and the integers $\sigma_i$ are called {\it structural indices}, or also {\it partial multiplicities}\footnote{Variations of this nomenclature can be found in the literature. For example, sometimes partial multiplicities are required to be nonzero integers (thus excluding zero indices), or even positive integers (in this case one distinguishes between the partial multiplicities as a zero, corresponding to the positive $\sigma_i$, and the partial multiplicities as a pole, corresponding to $-\sigma_j$ for the negative $\sigma_j$). In this paper, it is convenient to define the partial multiplicities as the exponents in the local Smith-McMillan form, thus admitting zero or negative integers.}, of $R(\lambda)$ at the point
$\lambda_0$. Moreover, they are ordered non-increasingly~: $\sigma_1 \geq \ldots \geq \sigma_r$. Observe that
negative indices are associated with poles of $R(\la)$, positive indices are associated with zeros of $R(\la)$, while zero indices are not associated to either. 

Other important sets of indices of a general $m\times n$ rational matrix 
$R(\la)$, called minimal indices, are related to its right null space $\ker R(\la)$ and left null space $\ker R(\la)^T$,
which are rational vector spaces over the field $\F(\la)$ of rational functions in $\la$.
\begin{definition} 
A matrix polynomial $N(\la) \in \F[\lambda]^{n\times p}$ of normal rank $p$
is called a minimal polynomial basis if the sum of the degrees of its columns, 
called the order of the basis, is minimal among all bases of 
span $N(\la)$. Its ordered column degrees are called the minimal indices of the basis.
\end{definition}
It was shown in \cite{For75} that the ordered list of indices is independent of 
the choice of minimal basis of the space. If we define the right null space $\ker R(\la)$ and left null space 
$\ker R(\la)^T$ of a $m\times n$  rational matrix $R(\la)$ of normal rank $r$ as 
the vector spaces of rational vectors $x(\la)$ and $y(\la)$ annihilated by $R(\la)$
$$  \ker R(\la):=\{ x(\la) \; | \; R(\la)x(\la)=0\}, \qquad    \ker R(\la)^T :=\{ y(\la) \;| \; y^\T(\la) R(\la)=0\}
$$
then the minimal indices of any minimal polynomial basis for these spaces, are called the right and left  
minimal indices of $R(\la)$. The respective dimensions of $ \ker R(\la)$ and $\ker R(\la)^T$ are $n-r$ and 
$m-r$ and the respective 
indices are denoted by
$$\{\epsilon_1,\ldots,\epsilon_{n-r}\}, \quad \{\eta_1,\ldots,\eta_{m-r}\}.
$$
It was also shown in \cite{For75} that for any minimal basis $N(\la)$, the constant matrix 
$N(\la_0)$ has full column rank for all $\la_0\in\F$ and the highest column degree matrix
of $N(\la)$ also has full column rank.

\subsection{Root polynomials} \label{Sec:zerodirections}

In \cite{DopN21,GLR82} a rigorous link was established between partial multiplicities and the so-called {\it root polynomials}
but in those references that concept was analyzed  for a polynomial matrix $P(\la)$ only. Root polynomials had appeared as technical tools for proving other results \cite{GLR82, N11}, and in \cite{DopN21} it was argued that they can play an important role in the theory of polynomial matrices. 
Below we give a definition which is equivalent to that given in \cite{DopN21}, but rephrased in a way that is more useful for deriving the results given in this paper. The matrix $P(\la)\in \F[\la]^{m\times n}$ is assumed to be polynomial in the next definition and theorem, proven in \cite{DopN21}.
\begin{definition} \label{def:rootpolynomial} Let $N(\la)$ be a right minimal basis of $P(\la)$ and $\la_0 \in \F$; then
\begin{itemize} \item $x(\la) \in \F[\la]^n$ is a root polynomial at $\la_0$ of order $k>0$ if $P(\la)x(\la)=(\la-\la_0)^k v(\la)$, $v(\la_0)\neq 0$, and $\left[ N(\la_0) \; x(\la_0)\right]$ has full column rank;
 \item $\{x_1(\la),...,x_s(\la)\}$ is a set of $\la_0$-independent root polynomials  of orders $\{k_1,...,k_s\}$ if
they are root polynomials at $\la_0$ of orders  $\{k_1,...,k_s\}$ and $\left[ N(\la_0) \; x_1(\la_0) \ldots x_s(\la_0) \right]$ has full column rank;
\item a $\la_0$-independent set is complete if there does not exist any $\la_0$-independent set of larger cardinality;
\item such a complete set is ordered if $k_1\ge \ldots \ge k_s > 0$;
\item such a complete ordered set is maximal if there is no root polynomial $\tilde x(\la)$ of order $k>k_j$ at 
$\la_0$ such that $\left[ N(\la_0) \; x_1(\la_0) \ldots x_{j-1}(\la_0) \; \tilde x(\la_0) \right]$ has full column rank for any $j$.
\end{itemize}
\end{definition}

The following result \cite[Theorem 4.1.3 and Theorem 4.2]{DopN21} explains the importance of these maximal sets.

\begin{theorem}\label{thm:frovan}
Let the nonzero partial multiplicities at $\lambda_0$ of $P(\la)$ be $0<\sigma_s \leq \dots \leq \sigma_1$, and suppose that  $x_1(\la),\dots,x_s(\la)$ are a complete set of root polynomials at $\lambda_0$ for $P(\la)$. Then, the following are equivalent:
\begin{enumerate}
\item $x_1(\la),\dots,x_s(\la)$ are a maximal set of root polynomials at $\lambda_0$ for $P(\la)$;
\item the orders of such a set are precisely $\sigma_1,\dots,\sigma_s$;
\item the sum of the orders of such a set is precisely $\sum_{i=1}^s \sigma_i$.
\end{enumerate}
\end{theorem}

In the next section, we show how one can extend the theory of root polynomials to rational matrix-valued functions, even when the point $\la_0$ is a pole. (If $\la_0$ is a zero, but not a pole, the task is trivially achieved by multiplying the rational matrix by its least common denominator, as the latter cannot have a root at $\la_0$ and, hence, it is easy to see that this operation preserves both partial multiplicities at $\la_0$ and maximal sets of root polynomials at $\la_0$.)

\section{Root vectors of rational matrices} \label{Sec:rootvectors}
In this section we build a theory of eigenvalues, (right) root vectors and (right) eigenvectors\footnote{The dual concepts of left root vectors and eigenvectors are, of course, trivially obtainable by working with $R(\la)^T$.} of any rational matrix $R(\la)$ over $\F(\la)$, without any need to make any of the assumptions that are sometimes made in the literature, e.g., that $R(\la)$ has full column rank or that the point of interest is not a pole. To this goal, we first need to recall what a local discrete valuation is.
\begin{definition}
Given $\la_0 \in \F$, the \emph{local discrete valuation $\kappa$} associated with the point $\la_0$ is the function
\[ \kappa : \F(\la) \rightarrow \Z \cup \{ \infty \} , \qquad 0\neq r(\la) = (\la-\la_0)^k \frac{p(\la)}{q(\la)} \mapsto \kappa[r(\la)]= k,  \quad 0 \mapsto \kappa[0]= \infty,\]
where, in the case of a nonzero rational function $r(\la)$, $p(\la)$ and $q(\la)$ are polynomials both coprime with $(\la-\la_0)$, and $k \in \Z$.
\end{definition} 
Note that, while there is a distinct local discrete valuation for each $\la_0$, in this paper we fix $\la_0$ and work locally. Thus, for simplicity of notation, we omit the dependence of the function $\kappa$ on the point $\la_0$. It is easy to check that the local discrete valuation at $\la_0$ satisfies the three properties
\begin{enumerate}
\item[(i)] $\kappa[a(\la)b(\la)] = \kappa[a(\la)] + \kappa[b(\la)]$, 
\item[(ii)] $\kappa[a(\la)+b(\la)] \geq \min \{ \kappa[a(\la)],\kappa[b(\la)] \}$, and 
\item[(iii)] $\kappa[a(\la)]=\infty \Leftrightarrow a(\la)=0$.
\end{enumerate}
Local discrete valuations can be extended to matrices by taking the elementwise minimum, and it retains certain (somewhat weaker) properties.
\begin{definition} If $A(\la) \in \F(\la)^{m \times n}$ and $\la_0 \in \F$, then the local discrete valuation at $\la_0$ of $A(\la)$ is $\kappa[A(\la)]=\min_{1 \leq i \leq m, 1 \leq j \leq n} \kappa[A_{ij}(\la)]$.
\end{definition}
\begin{proposition}\label{prop:ertiesofldv} For any pair of rational matrices $A(\la),B(\la)$ and any point $\la_0 \in \F$ with associated local discrete valuation $\kappa$, the following holds~:
\begin{enumerate}
\item if the product $P(\la)=A(\la)B(\la)$ is defined, then $\kappa[P(\la)] \geq \kappa[A(\la)]+\kappa [B(\la)]$;
\item  if the sum $S(\la)=A(\la)+B(\la)$ is defined, then $\kappa[S(\la)] \geq \min \{\kappa[A(\la)],\kappa[B(\la)]\}$;
\item  $\kappa[A(\la)]=\infty \Leftrightarrow A(\la)=0$;
\item  if $A(\la)$ is square and $\kappa[A(\la)]=\kappa[\det(A(\la))]=0$, then $A(\la)$ is invertible over $\F(\la)$ and \linebreak 
$\kappa[A(\la)^{-1}]=\kappa[\det(A(\la)^{-1})]=0$;
\item  if $A(\la)$ is square and satisfies $\kappa[A(\la)]=\kappa[\det(A(\la))]=0$, then \\ (a) if  the product $P(\la)=A(\la)B(\la)$ is defined, then $\kappa[P(\la)] =  \kappa[B(\la)]$ \\ (b) if the product $Q(\la)=B(\la)A(\la)$ is defined, then $\kappa[Q(\la)]=\kappa[B(\la)]$.  
\item if $S(\la)$ and $L(\la)$ are, respectively, the Smith-McMillan form of $A(\la)$ and the local (at $\la_0)$ Smith-McMillan form of $A(\la)$, then $\kappa[A(\la)]=\kappa[S(\la)]=\kappa[L(\la)]$;
\item if $A(\la)$ is partitioned into the two block submatrices $B(\la),C(\la)$ then $\kappa[A(\la)]=$\\ $\min \{\kappa[B(\la)], \kappa[C(\la)]\}$.
\end{enumerate}
\end{proposition}
\begin{proof}
\begin{enumerate}
\item For some indices $i,j$, 
\[ \kappa[P(\la)]=\kappa[P_{ij}(\la)] = \kappa \left[ \sum_k A_{ik}(\la)B_{kj}(\la) \right]\] \[\geq \min_k \kappa[A_{ik}(\la) ]+ \min_k \kappa [B_{kj}(\la)] \geq \kappa[A(\la)] + \kappa[B(\la)]. \]
\item \[ \kappa[A(\la) + B(\la)] = \min_{i,j} \kappa[A_{ij}(\la)+B_{ij}(\la)] \]
\[\geq \min_{i,j} \min \{ \kappa[A_{ij}(\la)] , \kappa[B_{ij}(\la)]\} \geq \min \{ \kappa[A(\la)],\kappa[B(\la)] \}.\]
\item Obvious.
\item If $A(\la)$ were not invertible, then it would follow that $\kappa[\det(A(\la))]=\infty$. Moreover, \[0=\kappa[\det(A(\la))]+\kappa[ \det(A(\la)^{-1})] \Rightarrow \kappa[\det(A(\la)^{-1})]=0\]
and
\[ 0 = \kappa[I] \geq \kappa[A(\la)] + \kappa[A(\la)^{-1}] \Rightarrow \kappa[A(\la)^{-1}] \leq 0. \]
Suppose now that $\kappa[A(\la)^{-1}]<0$. Then, since each entry of the inverse of $A(\la)$ is the ratio of the corresponding entry of the adjoint of $A(\la)$ (which must have nonnegative local discrete valuation being a polynomial in scalars with nonnegative local discrete valuation) and the determinant of $A(\la)$, we conclude that $\kappa[\det(A(\la))]\geq -\kappa[A(\la)^{-1}]>0$, contradicting 
$\kappa[\det(A(\la))]=0$.
\item We only prove part (a) as part (b) is similar. Under the assumptions, by the previous item $A(\la)$ is invertible and $\kappa[A(\la)^{-1}]=0$. Hence, by item 1,
\[ P(\la)=A(\la)B(\la) \Rightarrow \kappa[P(\la)] \geq \kappa[B(\la)]\]
and
\[ A(\la)^{-1}P(\la)=B(\la) \Rightarrow \kappa[B(\la)] \geq  \kappa [P(\la)],\]
yielding $\kappa[P(\la)] =\kappa[B(\la)]$. 
\item It suffices to prove $\kappa[S(\la)]=\kappa[A(\la)]$, as it is clear that $\kappa[S(\la)]=\kappa[L(\la)]$. We have that $U(\la) A(\la)=S(\la) V(\la)$ for some unimodular polynomial matrices $U(\la),V(\la)$. Clearly, being unimodular implies that $\kappa[\det( U(\la)) ]=\kappa[\det( V(\la))]=\kappa[U(\la)]=\kappa[V(\la)]=0$. Hence, by the previous item,
$$ \kappa[A(\la)]=\kappa[U(\la)A(\la)]=\kappa[S(\la)V(\la)]=\kappa[S(\la)].$$
\item Immediate by the definition of local discrete valuation of a matrix.
\end{enumerate}
\end{proof}

Let $R(\la) \in \F(\la)^{m \times n}$ and suppose that $N(\la)$ is a minimal basis for its right kernel. Given $\la_0 \in \F$, define as usual \cite{DopN21,N11} $\ker_{\la_0} R(\la) = \mathrm{span} N(\la_0)$. Then,
\begin{definition}\label{def:rootvec}
The rational vector $x(\la)$ is a \emph{root vector} at $\la_0$ of order $k>0$ for $R(\la)$ if:
\begin{enumerate}
\item $\kappa[x(\la)]=0$ and $x(\la_0) \not\in \ker_{\la_0} R(\la)$;
\item $\kappa[R(\la)x(\la)]=k>0$.
\end{enumerate}
Moreover, we say that $\la_0$ is an \emph{eigenvalue} (or zero) of $R(\la)$ if there exists a root vector at $\la_0$ for $R(\la)$.
\end{definition}

\begin{definition}
Let $N(\la)$ be a right minimal basis for $R(\la)$, then
\begin{itemize}
\item The set of root vectors $\{ x_i(\la) \}_{i=1}^s$ at $\la_0$, of orders $k_1,\dots,k_s$,  is $\la_0$-independent if the matrix
$\begin{bmatrix}
N(\la_0) & x_1(\la_0) & \dots & x_s(\la_0)
\end{bmatrix}$ has full column rank.
\item A $\la_0$-independent set of root vectors $\{x_i(\la) \}_{i=1}^s$ 
 is complete if there does not exist any $\la_0$-independent set $\{x_i(\la)\}_{i=1}^{\sigma}$  
 such that $\sigma>s$.
\item Such a complete set is ordered if $k_1 \geq \dots \geq k_s > 0$ 
\item Such a complete ordered set is maximal if, for any $j$, there is no root vector $\tilde{x}(\la)$ of order $k > k_j$ such that $$\begin{bmatrix}
N(\la_0) & x_1(\la_0) & \dots & x_{j-1}(\la_0) & \tilde{x}(\la_0)
\end{bmatrix}$$ has full column rank. 
\end{itemize} 
\end{definition}

At this point, we can state and prove Proposition \ref{prop:APB} below. Obviously, $\kappa[A(\la)]\ge 0$ implies that $\la_0$ is not a pole of $A(\la)$, and according to item 4 of Proposition \ref{prop:ertiesofldv}, $\kappa[A(\la)]=\kappa[\det(A(\la))]=0$ implies that $A(\la)$ is regular and that $\la_0$ is not an eigenvalue of $A(\la)$. Hence, Proposition \ref{prop:APB} describes the behaviour of root vectors when a rational matrix is multiplied (from both left and right) by square rational matrices whose local Smith form at $\la_0$ is trivial.

\begin{proposition}\label{prop:APB}
Let $R(\la),Q(\la) \in \F(\la)^{m \times n}$ satisfy $Q(\la)=A(\la)R(\la)B(\la)$ for some $A(\la) \in \F(\la)^{m \times m}$ and $B(\la) \in \F(\la)^{n \times n}$ such that $\kappa[A(\la)]=\kappa[B(\la)]=\kappa[\det(A(\la))]=\kappa[\det(B(\la))]=0$, where $\kappa$ is the local discrete valuation associated with $\la_0 \in \F$. Then: \begin{enumerate}
\item If $r(\la)$ is a root vector for $Q(\la)$ at $\la_0$ of order $k>0$, then $B(\la)r(\la)$ is a root vector for $R(\la)$ at $\la_0$ of the same order;
\item If $q(\la)$ is a root vector for $R(\la)$ at $\la_0$ of order $k>0$, then $B(\la)^{-1}q(\la)$ is a root vector for $Q(\la)$ at $\la_0$ of the same order;
\end{enumerate}
\end{proposition}
\begin{proof}
It suffices to prove item 1, since item 2 then follows from $R(\la)=A(\la)^{-1} Q(\la) B(\la)^{-1}$. By item 5 in Proposition \ref{prop:ertiesofldv}, $\kappa[r(\la)]=0 \Rightarrow \kappa[B(\la)r(\la)]=0$. By the same result, 
\[ k = \kappa [Q(\la) r(\la) ] = \kappa[ A(\la) R(\la) B(\la) r(\la)] = \kappa[R(\la)B(\la)r(\la)].\]
Now suppose for a contradiction that $B(\la_0)r(\la_0) \in \ker_{\la_0} R(\la)$. Then, by \cite[Lemma 2.9]{DopN21} there is a polynomial vector $v(\la)$ such that $R(\la)v(\la)=0$ and $v(\la_0)=B(\la_0)r(\la_0)$. (Although in \cite{DopN21} this is proved in the case where $R(\la)$ is a polynomial matrix, its proof is only based on the existence of a minimal basis for $\ker R(\la)$ which is still true when $R(\la)$ is generally rational.) This implies
\[ Q(\la) B(\la)^{-1} v(\la) = 0 \Rightarrow B(\la_0)^{-1} v(\la_0)=r(\la_0) \in \ker_{\la_0} Q(\la),\]
contradicting the assumption.
\end{proof}

We can then proceed similarly to \cite{DopN21}, and show the next theorem whose proof is omitted since it is completely analogous to the proof of \cite[Theorem 3.4]{DopN21}.
\begin{theorem}\label{thm:PtoAPB}
Let $R(\la),Q(\la) \in \F(\la)^{m \times n}$ satisfy $Q(\la)=A(\la)R(\la)B(\la)$ for some $A(\la) \in \F(\la)^{m \times m}$ and $B(\la) \in \F(\la)^{n \times n}$ such that $\kappa[A(\la)]=\kappa[B(\la)]=\kappa[\det(A(\la))]=\kappa[\det(B(\la))]=0$, where $\kappa$ is the local discrete valuation associated with $\la_0 \in \F$. Then: \begin{enumerate}
\item If $x_1(\la), \dots, x_s(\la)$  are a maximal (resp. complete, $\la_0$-independent) set for $Q(\la)$ with orders $k_1 \geq \dots \geq k_s > 0$, then $B(\la) x_1(\la),\dots,B(\la)x_s(\la)$ are a maximal (resp. complete, $\la_0$-independent) set for $R(\la)$ with the same orders;
\item If $x_1(\la), \dots, x_s(\la)$  are a maximal (resp. complete, $\la_0$-independent) set for $R(\la)$ with orders $k_1 \geq \dots \geq k_s > 0$, then $B(\la)^{-1} x_1(\la),\dots,B(\la)^{-1}x_s(\la)$ are a maximal (resp. complete, $\la_0$-independent) set for $Q(\la)$ with the same orders.
\end{enumerate}
\end{theorem}
At this point we observe that the set of all matrices $M(\la)$  satisfying the assumptions $\kappa[M(\la)]=\kappa[\det(M(\la))]=0$ (such as $A(\la)$ and $B(\la)$ in Theorem \ref{thm:PtoAPB}) are precisely the unimodular matrices over the local ring \begin{equation}\label{eq:localring}
{\cal R}:=\left\{ \frac{a(\la)}{b(\la)} \in \F(\la) : \gcd(a(\la),b(\la))=1 \ \mathrm{and} \ b(\la_0) \neq 0 \right\}.
\end{equation}

 It is clear that ${\cal R}$ is a principal ideal domain, and that the local Smith-McMillan form at $\la_0$ is obtained precisely by unimodular equivalence over ${\cal R}$. Hence, one can reduce the study of maximal sets to the case of rational matrices in local Smith-McMillan form. In particular, the proofs of \cite[Proposition 2.16 and Theorems 3.5, 4.1 and 4.2]{DopN21} generalize to the rational case, as we show below.
\begin{definition} \label{kermu}
 Let $R(\la) \in \F(\la)^{m \times n}$. 
 Then, 
\[ \ker R(\la_0) = \{ x_0 \in \F^n : \exists x(\la) \in \F(\la)^n \ s.t. \ \kappa[R(\la)x(\la)]>0 \ \mathrm{and} \ x(\la_0)=x_0 \} \]
where $\kappa$ is the local discrete valuation at $\la_0$.\end{definition}

Note that the condition $x(\la_0)=x_0$ in Definition \ref{kermu} implies in particular that $\la_0$ is not a pole of $x(\la)$, or equivalently that $\kappa[x(\la)]\geq 0$ (and $=0$ if and only if $x_0 \neq 0$). It is not hard to verify that $\ker R(\la_0)$, as defined above, is a subspace of $\F^n$.

\begin{remark}
In the case where $\la_0$ is not a pole of $R(\la)$, i.e., $R(\la_0) \in \F^{m \times n}$, Definition \ref{kermu} is equivalent to 
$$\ker R(\la_0) =\{ x_0 \in \F^n :  R(\la_0)x_0=0\}.$$ 
This is a more familiar expression but, in general, it does not make sense because $R(\la_0)$ may not be well defined. In practice, usually the most relevant case is when $\F=\C$ or some subfield of $\C$: in this case, Definition \ref{kermu} is equivalent to 
$$\ker R(\la_0) =\{ x_0 \in \F^n : \exists x(\la) \in \F(\la)^n \ s.t. \lim_{\lambda \rightarrow \la_0}R(\la)x(\la)=0 \ \mathrm{and} \ x(\la_0)=x_0\}$$ 
since $\kappa[R(\la)x(\la)] > 0 \Leftrightarrow \lim_{\lambda \rightarrow \la_0}R(\la)x(\la)=0$. Related limit-based definitions have appeared elsewhere, at least under the simplifying assumption that $R(\la)$ has full column rank, to characterize eigenvectors of rational matrices with the additional condition $x_0 \neq 0$ (clearly such a nonzero $x_0$ exists only if $\ker R(\la_0)$ is not the trivial subspace): see for instance \cite[p. 315]{MITnotes}. We will give a more general characterization of eigenvectors in Definition \ref{def:evec}. We also point out that Definition \ref{kermu} is more general than the limit-based definitions such as the one in \cite[p. 315]{MITnotes} and, being purely algebraic, is valid also when $\F$ is not a topological space; this choice is coherent with the literature on polynomial matrices that have often been studied over any field.
\end{remark}

\begin{theorem}\label{span}
Let $R(\la) \in \F(\la)^{m \times n}$ with minimal basis $N(\la) \in \F[\la]^{n \times p}$. Let $\{x_i(\la)\}_{i=1}^s$ be root vectors at $\la_0$ for $R(\la)$; then, they are a complete set if and only if the columns of the matrix
\[ M:=\begin{bmatrix} N(\la_0) & x_1(\la_0) & \dots & x_s(\la_0) \end{bmatrix}\]
are a basis for $\ker R(\la_0)$. \end{theorem}
\begin{proof} By definition, $\kappa[R(\la)N(\la)]=\infty$ and $\kappa[R(\la)x_i(\la)]>0$. Hence, the columns of $M$ are in $\ker R(\la_0)$. Moreover, if the set $\{x_i(\la)\}_{i=1}^s$ is complete, $M$ has full column rank as any complete set is $\la_0$-independent. Suppose the columns of $M$ are not a basis for $\ker R(\la_0)$: then they can be completed to a basis, say, by adding some vectors $\{ v_i \}_{i=1}^\ell$. But then there are rational vectors $w_i(\la)$ satisfying both $w_i(\la_0)=v_i$ and $\kappa[ R(\la) w_i(\la) ] > 0$. This implies that $\{x_i(\la)\}_{i=1}^s \cup \{w_i(\la)\}_{i=1}^\ell$ is still a $\la_0$-independent set. Hence, the original set was not complete. Conversely, if the set was not complete, then the matrix $M$ cannot be a basis for $\ker R(\la_0)$ because then by starting from a complete set we can find a full column rank matrix $\hat M$, having strictly more columns than $M$, and such that the columns of $\hat M$ belong to $\ker R(\la_0)$.
\end{proof}

We observe that Theorem \ref{span} in particular implies that the cardinality of the subset of root vectors in a complete set is equal to $s=\dim \ker R(\la_0)-\dim \ker{_{\la_0}} R(\la)$.

\begin{theorem}\label{span2}
Every complete set of root vectors has the same cardinality. And in every maximal set, the orders of the root vectors are the same.
\end{theorem}
\begin{proof}
The first statement is obvious by Theorem \ref{span}. The second statement follows from 
the proof of \cite[Theorem 4.1.3]{DopN21}, which shows that root polynomials in a maximal set must have the same orders when $R(\la)$ is polynomial; however, it can be applied verbatim to the rational case. (More details are in the proof of Theorem \ref{span3} below.) 
\end{proof}

\medskip

We can now formulate the major result of this section. It extends \cite[Theorem 4.1]{DopN21} to the rational case and has essentially the same proof; given its importance, we include a sketch of proof below.
\begin{theorem}\label{span3}
Let $R(\la) \in \F(\la)^{m \times n}$ and suppose that the positive partial multiplicities of $\la_0 \in \F$ are $\sigma_1 \geq \dots \geq \sigma_s > 0.$ If $R(\la)$ has a complete set of root vectors at $\la_0$ $x_1(\la),\dots,x_s(\la)$, having orders $ k_1 \geq \dots \geq  k_s$, then $\sigma_i \geq k_i$ for all $i$, and such a set is maximal if and only if $k_i = \sigma_i$ for all $i=1,\dots,s$.
\end{theorem}
\begin{proof}
Following \cite{DopN21}, let $\{x_i(\la)\}_{i=1}^s$ be a complete set and $\{v_i(\la)\}_{i=1}^s$  be a maximal set, and suppose that these sets have orders $k_1 \geq \dots \geq k_s$ and $\tau_1 \geq \dots \geq \tau_s$, respectively. Suppose for a contradiction that for some $j \geq 0$ it holds $k_i=\tau_i$ for $i\leq j$ but $k_{j+1} > \tau_{j+1}$. Let $\{r_i(\la)\}_{i=1}^p$ be a minimal basis for $R(\la)$, then by the definition of maximal set the matrix
\[ \begin{bmatrix}
r_1(\la_0) & \dots & r_p(\la_0) & v_1(\la_0) & \dots & v_j(\la_0) & x_k(\la_0)
\end{bmatrix}\]
cannot possibly have full column rank, for all $k=1,\dots,j+1$. This implies that \[ r_1(\la_0),\dots, r_p(\la_0),x_1(\la_0),\dots,x_{j+1}(\la_0) \in \mathrm{span}\{ r_1(\la_0),\dots,r_p(\la_0),v_1(\la_0),\dots,v_j(\la_0)\};\]
but the latter claim exhibits $p+j+1$ linearly independent (by completeness and hence $\la_0$-independence of $\{x_i(\la)\}_{i=1}^s$) vectors all belonging to the same $(p+j)$-dimensional vector space, which is absurd. This shows that the orders of a complete set are elementwise bounded above by the orders of a maximal set. Moreover, if we further assume that also $\{x_i(\la)\}_{i=1}^s$ is maximal, then by a symmetric argument the case where $k_{j+1} < \tau_{j+1}$ can also be excluded analogously: hence, any two maximal sets have the same orders.

To conclude the proof, it suffices to exhibit a particular maximal set and check that its orders are the positive partial multiplicities at $\la_0$. To this goal, let $S(\la)=U(\la)R(\la)V(\la)$ be the local Smith-McMillan form at $\la_0$ of $R(\la)$ for some $U(\la),V(\la)$ invertible over the local ring \eqref{eq:localring}, and suppose that the integers $m,s$ are such the positive partial powers of $(\la-\la_0)$ appear precisely in positions $1 \leq m, \dots, m+s-1$. It is a simple exercise to prove that $\{e_{i+1-m}\}_{i=1}^s$ is a maximal set of root vectors for $S(\la)$ (see also \cite[Theorem 3.1]{DopN21} and note that the proof extends to the rational case) and that their orders are precisely the positive partial multiplicities at $\la_0$. By Theorem \ref{thm:PtoAPB}, $\{ V(\la) e_{i+1-m}\}_{i=1}^s$ is a maximal set of root vectors for $R(\la)$ with the same orders.
\end{proof}

In other words, also in the case of rational matrices, maximal sets of root vectors yield the full information about positive partial multiplicities.

\begin{remark} We point out that the root vectors have their entries in the ring \eqref{eq:localring}. They can therefore always be scaled with an appropriate polynomial factor to make them polynomial root vectors. These scaling factors do not affect any of the theorems derived in this section.
\end{remark}

From the theory of root vectors for rational matrices, one can define eigenvectors by generalizing the approach of \cite[Definition 2.18]{DopN21}. We restrict ourselves to right eigenvectors; left eigenvectors can be defined as right eigenvectors of $R(\la)^T$. Namely, eigenvectors are nonzero vectors in the quotient space of $\ker R(\la_0)/\ker_{\la_0} R(\la)$ (which has a vector space structure). Observe that, if $R(\la)$ has full column rank, then $\ker_{\la_0} R(\la)=\{0\}$ and hence $\ker R(\la_0)/\ker_{\la_0}R(\la)=\ker R(\la_0)$.

\begin{definition}\label{def:evec} An equivalence class $[v] \in \ker R(\la_0)/\ker_{\la_0} R(\la)$ is called a (right) eigenvector of $R(\la)$ associated with $\la_0$ if $[v] \neq [0].$
\end{definition} 
Definition \ref{def:evec} implies that eigenvectors associated with $\la_0$ exist if and only if $\la_0$ is an eigenvalue, as by Definition \ref{def:rootvec} and Definition \ref{kermu} this is equivalent to the inclusion $\ker_{\la_0} R(\la) \subset \ker R(\la_0)$ being strict; if, on the contrary, $\la_0$ is not an eigenvalue then $\ker_{\la_0}R(\la)=\ker R(\la_0)$, and hence the quotient space is the trivial vector space $\{[0]\}$ which does not contain nonzero equivalence classes (and thus it cannot contain eigenvectors). Moreover, by Theorem \ref{span}, it is clear that $[v]$ is an eigenvector if and only if $$[0]\neq [v] \in \mathrm{span}( [x_1(\la_0)], \dots, [x_s(\la_0)] )$$ where $\{x_i(\la)\}_{i=1}^s$ is any complete set of root vectors at $\la_0$ for $R(\la)$. In other words, $[v]$ is an eigenvector if and only if $v=\sum_{i=1}^s x_i(\la_0) c_i + w$ for some $c_i \in \F$ such that $\sum_{i=1}^s x_i(\la_0) c_i\neq 0$, and some $w \in \ker_{\la_0} R(\la)$. When $\F$ is a subfield of $\C$ and $\la_0$ is simple, one can uniquely (up to a nonzero scalar) fix a representative $u$ of the equivalence class $[v]$ by asking that $u \in \ker_{\la_0} R(\la)^\perp$; see \cite{LN20}. Finally, we emphasize again that if $R(\la)$ has full column rank then $\ker_{\la_0} R(\la)$ is trivial, and therefore Definition \ref{def:evec} collapses to an eigenvector being equal to $v=x(\la_0)$ where $x(\la)$ is a root vector at $\la_0$ for $R(\la)$ (note that by definition root vectors do not have poles at $\la_0$ and are nonzero when evaluated at $\la_0$). If, furthermore, $\la_0$ is not a pole, then we recover the classical definition: a nonzero vector $v$ such that $R(\la_0)v=0$.

\subsection{Root vectors at infinity}\label{sec:infinity}
We conclude this section by noting that, as in the polynomial case \cite{DopN21}, it is possible to define root vectors at $\infty$ for rational matrices. There are two possible approaches for this. One, natural for polynomial eigenvalue problems, is via the reversal and was described in \cite[Sec. 6]{DopN21}. A second one, more natural in the context of control theory as it reveals the Smith-McMillan positive degrees at infinity,  is via the local discrete valuation at infinity of a rational function \cite{For75,vardoulakis}. The latter is defined for scalar rational functions as $\kappa[0]=\infty$ for the zero function, and  $\kappa[p(\la)/q(\la)]=\deg q(\la)-\deg p(\la)$ if $p(\la)/q(\la) \neq 0$.  It is a simple exercise to verify that this is indeed a local discrete valuation in the sense of satisfying the same three properties as the valuation at a finite point. We can, again, extend the function $\kappa$ to rational matrices by taking the elementwise minimum of its value on the matrix entries. A root vector $x(\la)$ of order $k$ at $\infty$ for $R(\la)$ is then a rational vector such that 
\begin{enumerate}
    \item  $\kappa[x(\la)]=0$, that is, $x(\la)$ is a proper rational vector with at least one biproper entry;
     \item $x(\infty)\not \in \ker_\infty R(\la)$ where (2a) $x(\infty)$ is defined as the zeroth-order coefficient $x_0$ in the expansion of $x(\la)=\sum_{j=0}^\infty x_j \la^{-j}$ (which exists because $x(\la)$ is proper) and (2b) $\ker_\infty R(\la):=\ker_0 R(\frac{1}{\la})$ or, equivalently \cite{N11}, the vector space spanned by the high order coefficients matrix of a minimal basis for $\ker R(\la)$ ;
     \item   $\kappa[R(\la)x(\la)] = k >0$, that is, $R(\la)x(\la)$ is a strictly proper rational vector with local discrete valuation at $\infty$ equal to $k$.
\end{enumerate}
With this definition, $x(\la)$ is a root vector of order $k$ at $\infty$ for $R(\la)$ if and only if $x(\la^{-1})$ is a root vector of order $k$ at $0$ for $R(\la^{-1})$. Moreover, $\infty$-independent, complete and maximal sets for $R(\la)$ can then be defined similarly to the finite case, and correspond to $0$-independent, complete and maximal sets for $R(\la^{-1})$.  Hence, the properties of root vectors that we have established for finite points also carry over to $\infty$, provided that the partial multiplicities at $\infty$ are defined \cite[Sec. 5]{amzinf} as the ones obtained via the Smith form over the local ring of proper rational functions. For reasons of conciseness, we omit the details. We warn the reader that, as we illustrate below, this definition of partial multiplicities at infinity (common in the context of control theory) differs from the one usually given in the context of polynomial eigenvalue problems.

\begin{example}
It should be noted that the definition of a root vector at infinity proposed above does not yield the one given in \cite{DopN21} in the case of polynomial matrices. The reason is that the partial multiplicities at infinity are defined differently in the two approaches (via the reversal or via the Smith-McMillan approach).

For instance consider the regular polynomial matrix $P(\la)=\begin{bmatrix}1&\la\\
0&1\end{bmatrix}$. According to the reversal definition (and seeing $P(\la)$ as grade $1$), the partial multiplicities of $P(\la)$ at infinity are $0,2$; according to the Smith-McMillan definition, however, they are $-1,1$. This difference is reflected in the two different definitions of root polynomial (or root vector) at infinity given in \cite{DopN21} or here. According to the definition in \cite{DopN21}, $r(\la)=\begin{bmatrix}-\la&1\end{bmatrix}^T$ is a root polynomial at $\infty$ of order $2$ for $P(\la)$ (and a maximal set) because it is the $1$-reversal of a root polynomial of order $2$ of the $2$-reversal of $P(\la)$. According to the definition proposed above, instead, $r(\la)$ is not even a root vector at $\infty$ for $P(\la)$ (its local discrete valuation at $\infty$ is $-1\neq0$). However, $\la^{-1}r(\la)$ is a root vector of order $1$ at $\infty$ since (1) $\kappa[\la^{-1}r(\la)]=0$ (2) $r(\infty)=\begin{bmatrix}-1&0\end{bmatrix}^T \not \in \ker_\infty P(\la)=\{0\}$ (3) $\kappa[P(\la)\la^{-1}r(\la)]=1$.
\end{example}

\section{The system matrix of a rational matrix}  \label{Sec:systemmatrix}

Any rational transfer function $R(\lambda)$ can be realized by a {\it
generalized state-space realization} $\{ A,B,C,D,E\}$  (see e.g.\ \cite{kailath}) such that
\begin{equation}\label{eq:wrongref}
R(\lambda) = C(\lambda E-A)^{-1}B+D.
\end{equation}
If $R(\lambda)$ happens to be proper, i.e. bounded at $\la=\infty$,
then $E$ can be chosen to be the identity and one refers to the system
$\{A,B,C,D\}$ as a {\it state-space realization} of $R(\lambda)$.
One often associates with the above realization \eqref{eq:wrongref} the pencil
\begin{equation}
S(\lambda) := \left[ \begin{array}{c|c} A-\lambda E & B \\
\hline C & D \end{array} \right]
\end{equation}
called the {\it system matrix} of the realization. Note that $R(\lambda)$
is the Schur complement of its system matrix $S(\la)$. Techniques for constructing a realization
\eqref{eq:wrongref} are described in e.g. \cite{vd81}. It is also shown there how
to obtain an {\it irreducible generalized state space realization} of the transfer function $R(\lambda)$. 
These are realizations \eqref{eq:wrongref} for which
\begin{equation} \label{irreduce}
\left[ \begin{array}{c} A-\lambda E  \\
\hline C \end{array} \right]
\quad \mathrm{and} \quad 
\left[ \begin{array}{c|c} A-\lambda E & B
\end{array} \right]
\end{equation}
have no finite zeros (see also \cite{vvk}).

\smallskip
More generally, one can consider a \emph{polynomial system matrix}
\begin{equation}\label{eq:psm}
 P(\la) = \begin{bmatrix}
-A(\la) & B(\la) \\
C(\la) & D(\la)
\end{bmatrix} \in \F[\la]^{(q+m) \times (q+n)} \end{equation}
such that $A(\la) \in \F[\la]^{q \times q}$ is regular and $R(\la)=D(\la)+C(\la)A(\la)^{-1}B(\la)$. The polynomial system matrix is said to be \emph{minimal} if the Smith normal form of $\begin{bmatrix}
-A(\la) & B(\la)
\end{bmatrix}$ is $\begin{bmatrix}
I&0
\end{bmatrix}$ and the Smith normal form of $\begin{bmatrix}
-A(\la)\\
C(\la)
\end{bmatrix}$ is $\begin{bmatrix}
I\\
0
\end{bmatrix}$. In the following section, we will investigate how root vectors of a rational matrix and its minimal polynomial system matrix relate.

\section{Root vectors of a minimal system matrix and those of its transfer function}  \label{Sec:rootsofSandR}

As was shown in \cite{Rosenbrock}, much of the structure of $R(\lambda)$
can be retrieved from that of $P(\lambda)$ as in \eqref{eq:psm}, when the latter is a minimal system matrix. It is shown there
that in that case
\begin{itemize}
\item the pole structure of $R(\lambda)$ at its finite 
poles is identical to the zero structure of $A(\lambda)$ at its finite zeros.
\item the zero structure of $R(\lambda)$ at its finite zeros
is identical to the zero structure of $P(\lambda)$ at its finite  zeros.
\item the left and right null-space structures of $R(\lambda)$ and
$P(\lambda)$ are the same.
\end{itemize}
In the above relations, {\it structure} refers to the index sets
of the finite elementary divisors and of the left and right
minimal bases. In the rest of this section we show there is also a relation
between the {\it root vectors} as well as the index sets. For simplicity, we focus on finite eigenvalues.

\medskip

\begin{theorem}\label{thm:extraction}
Let $R(\la)=D(\la)+C(\la)A(\la)^{-1}B(\la) \in \F(\la)^{m \times n}$ be a rational matrix and let $P(\la)$ as in \eqref{eq:psm}
be a minimal polynomial system matrix of $R(\la)$.
If $x(\la)$ is a root vector of order $k$ at $\la_0$ for $R(\la)$, then $y(\la)=\begin{bmatrix}
A(\la)^{-1}B(\la)x(\la)\\
x(\la)
\end{bmatrix}$ is a root vector of order $k$ at $\la_0$ for $P(\la)$.

Moreover, if $\{x_i(\la)\}_{i=1}^s$ are a maximal set of root vectors at $\la_0$ for $R(\la)$ and $y_i(\la)=\begin{bmatrix}
A(\la)^{-1}B(\la)x_i(\la)\\
x_i(\la)
\end{bmatrix}$ then $\{y_i(\la)\}_{i=1}^s$ are a maximal set of root vectors at $\la_0$ for $P(\la)$.
\end{theorem}
\begin{proof}
It is readily verified that
$$ P(\la) y(\la) = \begin{bmatrix}
0\\
R(\la)x(\la)
\end{bmatrix} \Rightarrow \kappa[P(\la)y(\la)]=k.$$
Since by assumption $\kappa[x(\la)]=0$, to show that $\kappa[y(\la)]=0$ it suffices to verify that $\kappa[A(\la)^{-1}B(\la)x(\la)] \geq 0$; this can be proved as in \cite[Proposition 2.5]{NNPQ22} (the result in \cite{NNPQ22} is stated for $\F=\C$ and proved using limits, but it is not difficult to extend it to any field and by using valuations). At this point, the proof of \cite[Theorem 3.8, part (a)]{admz} shows that, if $M(\la) \in \F[\la]^{n \times p}$ is a minimal basis for $\ker R(\la)$, then $N(\la)=\begin{bmatrix}
A(\la)^{-1}B(\la)M(\la)\\
M(\la)
\end{bmatrix}$ is a basis for $\ker P(\la)$ that has full column rank for all $\la \in \F$. As a consequence, $x(\la_0)\not\in \ker_{\la_0} R(\la) \Rightarrow y(\la_0) \not\in \ker_{\la_0} P(\la)$.

Suppose now that $\{x_i(\la)\}_{i=1}^s$ are a maximal set for $R(\la)$. Then, $\{y_i(\la) \}_{i=1}^s$ are $\la_0$-independent because
\[ p+s \geq \rank \begin{bmatrix}
N(\la_0) & y_1(\la_0) & \dots & y_s(\la_0) \end{bmatrix} \geq \rank \begin{bmatrix}
M(\la_0) & x_1(\la_0) & \dots & x_s(\la_0) \end{bmatrix} = p+s.\]
 That $\{y_i(\la)\}_{i=1}^s$ are also complete (resp. maximal) follows by Theorem \ref{span} (resp. Theorem \ref{span3}) and by Rosenbrock's theorem \cite{C72}.
\end{proof}

A natural question is whether the converse of Theorem \ref{thm:extraction} holds, i.e., if the bottom block of a root polynomial for $P(\la)$ is always a root vector for $R(\la)$. The answer is positive if we consider a point $\la_0$ which is not a pole of $R(\la)$.

\begin{proposition}\label{prop:conditionalconverse}
Let $R(\la)=D(\la)+C(\la)A(\la)^{-1}B(\la)$ be a rational matrix and let $P(\la)$ as in \eqref{eq:psm}
be a minimal polynomial system matrix of $R(\la)$. Suppose, moreover, that $\la_0$ is not a pole of $R(\la)$.
If $y(\la)$ is a root vector of order $k$ for $P(\la)$, then there exists a rational vector $v(\la)$ such that $\kappa[v(\la)]\geq 0$, $y(\la)=\begin{bmatrix}
A(\la)^{-1}B(\la)x(\la) + (\la-\la_0)^k v(\la)\\
x(\la)
\end{bmatrix}$, and $x(\la)$ is a root vector of order at least $k$ at $\la_0$ for $R(\la)$.

Moreover, if $\{y_i(\la)\}_{i=1}^s$ are a maximal set of root vectors at $\la_0$ for $P(\la)$ and for all $i$ we denote by $x_i(\la)$ the bottom block of $y_i(\la)$, then  $\{x_i(\la)\}_{i=1}^s$ are a maximal set of root vectors at $\la_0$ for $R(\la)$.
\end{proposition}
\begin{proof}
Since $P(\la)$ is minimal and $\la_0$ is not a pole of $R(\la)$, it follows by \cite{C72} that $\la_0$ is not an eigenvalue of $A(\la)$. Hence, $A(\la)$ is unimodular over the local ring \eqref{eq:localring} and so is 
$ \begin{bmatrix}
I & 0\\
C(\la)A(\la)^{-1}&I
\end{bmatrix}.$ Thus, by Proposition \ref{prop:APB}, $y(\la)$ is a root vector at $\la_0$ of order $k$ for 
$$ \begin{bmatrix}
I & 0\\
C(\la)A(\la)^{-1}&I
\end{bmatrix} P(\la) = \begin{bmatrix}
-A(\la)&B(\la)\\
0&R(\la)
\end{bmatrix}.$$

Decomposing $y(\la)=\begin{bmatrix}
A(\la)^{-1}B(\la)x(\la)+z(\la)\\
x(\la)
\end{bmatrix}$ it follows that
$$ \kappa \left[ \begin{bmatrix} -A(\la)z(\la) \\
R(\la)x(\la) \end{bmatrix} \right] = k.$$
Hence, $\kappa[R(\la)x(\la)]\geq k$. Moreover, by item 5 in Proposition \ref{prop:ertiesofldv},
$\kappa[z(\la)]=\kappa[-A(\la)z(\la)] \geq k$ and therefore we can write $z(\la)=(\la-\la_0)^k v(\la)$ as sought. Furthermore, it must be $\kappa[x(\la)]=0$, for if the valuation was positive then we would have that \[ \kappa[y(\la)] \geq \min \{ \kappa[x(\la)],\kappa[B(\la)x(\la)],\kappa[z(\la)]\} \geq \min \{ k, \kappa[x(\la)]\}>0. \]
The argument to show $y(\la_0)\not\in \ker_{\la_0} P(\la) \Rightarrow x(\la_0) \not\in \ker_{\la_0} R(\la_0)$ is similar to its analogue in the proof of Theorem \ref{thm:extraction}, except that this time we need to apply the proof of \cite[Theorem 3.6, part (a)]{admz}.

Suppose now that $\{y_i(\la)\}_{i=1}^s$ are a maximal set. We claim that $\{x_i(\la)\}_{i=1}^s$ are $\la_0$-independent: this suffices to prove that they are maximal, by Rosenbrock's Theorem \cite{C72}, Theorem \ref{span} and Theorem \ref{span3}. To prove the claim, suppose for a contradiction that it is false. Denote by $M(\la)$ a minimal basis of $\ker R(\la)$ and observe \cite{admz} that $N(\la) = \begin{bmatrix}
A(\la)^{-1}B(\la)M(\la)\\
M(\la)\end{bmatrix}$ is a basis for $\ker P(\la)$ that has full rank for $\la=\la_0$. Then, for some nonzero constant vector $c$, it holds $0=\begin{bmatrix}
M(\la_0) & x_1(\la_0) & \dots x_s(\la_0)
\end{bmatrix} c$
which implies \[\begin{bmatrix}A(\la_0)^{-1}B(\la_0)\\
I \end{bmatrix}\begin{bmatrix}
M(\la_0) & x_1(\la_0) & \dots x_s(\la_0)
\end{bmatrix} c = \begin{bmatrix}
N(\la_0) & y_1(\la_0) & \dots y_s(\la_0)
\end{bmatrix} c = 0.\]
\end{proof}

(A special case of) Proposition \ref{prop:conditionalconverse} shows that the problem of finding a maximal set of root vectors of a
transfer function at a point which is a zero, but not a pole, can be reduced to the corresponding problem of a pencil
\begin{equation}
S(\lambda) =
\left[ \begin{array}{c|c} A & B \\
\hline C & D \end{array} \right] - \lambda
\left[ \begin{array}{c|c} E & 0 \\
\hline 0 & 0 \end{array} \right]
\end{equation}
derived from an irreducible generalized state-space realization. When $\F\subseteq\C$, an algorithm to compute root vectors (in fact, root polynomials, which are of course a special case of root vectors satisfying the extra constraint of being polynomial) of a pencil was recently studied in \cite{NofV}.

If we remove the assumption that $\la_0$ is not a pole, then the statement of Proposition \ref{prop:conditionalconverse} may no longer be true as the following example illustrates.

\begin{example}\label{ex:bad}

$$ R(\la) = \begin{bmatrix}
1 & 0\\
\frac{1}{\la} & 1
\end{bmatrix}, \qquad P(\la) = \begin{bmatrix}
-\la & 1 &0 \\
0&1&0\\
1&0&1
\end{bmatrix}.$$
It is easy to verify that $P(\la)$ is a minimal polynomial system matrix of $R(\la)$. The vector $\begin{bmatrix}
-1&0&1
\end{bmatrix}^T$ is a root polynomial of order $1$ at $0$ for $P(\la)$. However, its bottom block $\begin{bmatrix}
0&1
\end{bmatrix}^T$ is not a root vector at $0$ for $R(\la)$.
\end{example}

To overcome this difficulty, we need a technique to remove the poles while keeping the original positive partial multiplicities, and simultaneously being able to recover maximal sets of root vectors. Let us start from a sufficient condition for root polynomial recovery even in the presence of poles.

\begin{lemma}\label{lem:sufficient}
Let $R(\la)=D(\la)+C(\la)A(\la)^{-1}B(\la)$ be a rational matrix and let $P(\la)$ as in \eqref{eq:psm}
be a minimal polynomial system matrix of $R(\la)$.
Suppose that $y(\la)=\begin{bmatrix}
u(\la)\\
x(\la)\end{bmatrix}$ is a root vector of order $k$ at $\la_0$ for $P(\la)$ satisfying the additional condition $A(\la)u(\la)=B(\la)x(\la)$. Then, $x(\la)$ is a root vector of order $k$ at $\la_0$ for $R(\la)$.

Moreover, suppose that $\{ y_i(\la) \}_{i=1}^s$ is a maximal set of root vectors at $\la_0$ for $P(\la)$, all satisfying the additional condition $A(\la)u_i(\la)=B(\la)x_i(\la)$ where $y_i(\la)=\begin{bmatrix}u_i(\la)\\
x_i(\la) \end{bmatrix}$. Then, $\{ x_i(\la) \}_{i=1}^s$ is a maximal set of root vectors at $\la_0$ for $R(\la)$.
\end{lemma}
\begin{proof}
By assumption, $u(\la)=A(\la)^{-1}B(\la)x(\la)$ and hence \[ P(\la)y(\la) = \begin{bmatrix}
0\\
R(\la)x(\la)\end{bmatrix} \Rightarrow \kappa[R(\la)x(\la)]=\kappa[P(\la)y(\la)]=k. \]
Suppose now for a contradiction that $x(\la_0) \in \ker_{\la_0} R(\la)$; then, if $M(\la)$ is a minimal basis for $\ker R(\la)$, $x(\la_0)=M(\la_0)c$ for some constant vector $c$. On the other hand, by the results in \cite{admz}, \[ N(\la) = \begin{bmatrix}
A(\la)^{-1}B(\la)M(\la)\\
M(\la) \end{bmatrix} \]
is a polynomial basis for $\ker P(\la)$ that has full column rank upon evaluation at $\la=\la_0$, and therefore we obtain the contradiction $y(\la_0)=N(\la_0)c \Rightarrow y(\la_0) \in \ker_{\la_0} P(\la)$. Finally, clearly $\kappa[x(\la)] \geq \kappa[y(\la)]=0$, but if $\kappa[x(\la)]>0$ then $x(\la_0)=0 \in \ker_{\la_0} R(\la)$, and thus $\kappa[x(\la)]=0$.

Now suppose that $\{y_i(\la)\}_{i=1}^s$ are a maximal set of root vectors at $\la_0$ for $P(\la)$ as in the statement. We claim that $\{x_i(\la)\}_{i=1}^s$ are $\la_0$-independent. Hence, by Rosenbrock's Theorem \cite{C72}, Theorem \ref{span} and Theorem \ref{span3}, we conclude that they are also maximal. It remains to prove the claim. To this goal, denote
\[ X(\la):=\begin{bmatrix}x_1(\la)&\dots&x_s(\la)\end{bmatrix}, \qquad Y(\la):=\begin{bmatrix}y_1(\la)&\dots&y_s(\la)\end{bmatrix}\]
and define $F(\la):=A(\la)^{-1}B(\la)M(\la)$, $U(\la):=A(\la)^{-1}B(\la)X(\la)$; as already observed, neither $F(\la)$ nor $U(\la)$ have poles at $\la_0$. Note also that, since $\begin{bmatrix}-A(\la)\\
C(\la)\end{bmatrix}$ has trivial Smith form, it has some polynomial left inverse $L(\la)$. Moreover, there is a matrix $W(\la)$ with the property that $W(\la_0)=0$ and such that
\[ \begin{bmatrix}0&0\\
0&W(\la)\end{bmatrix}=P(\la) \begin{bmatrix}N(\la)&Y(\la)\end{bmatrix}=\begin{bmatrix}-A(\la)\\
C(\la)\end{bmatrix} \begin{bmatrix}F(\la)&U(\la)\end{bmatrix} +\begin{bmatrix}B(\la)\\
D(\la)\end{bmatrix}\begin{bmatrix}M(\la)&X(\la)\end{bmatrix}\]
which implies
\[ \begin{bmatrix}F(\la_0)&U(\la_0)\end{bmatrix}=-L(\la_0)\begin{bmatrix}B(\la_0)\\
D(\la_0)\end{bmatrix}\begin{bmatrix}M(\la_0)&X(\la_0)\end{bmatrix}=:H \begin{bmatrix}M(\la_0)&X(\la_0)\end{bmatrix} \]
where $H$ is some constant matrix. Thus,
\[ \begin{bmatrix}N(\la_0)&Y(\la_0)\end{bmatrix}=\begin{bmatrix}F(\la_0)&U(\la_0)\\M(\la_0)&X(\la_0)\end{bmatrix}=\begin{bmatrix}H\\
I\end{bmatrix}\begin{bmatrix}M(\la_0)&X(\la_0)\end{bmatrix}.\] Suppose now for a contradiction that there exists a nonzero constant vector $c$ sastisfying $\begin{bmatrix}M(\la_0)&X(\la_0)\end{bmatrix}c=0$; then,
\[ \begin{bmatrix}N(\la_0)&Y(\la_0)\end{bmatrix}c = \begin{bmatrix}H\\
I\end{bmatrix} 0 = 0,  \]
contradicting the $\la_0$-independence of $\{y_i(\la)\}_{i=1}^s$.
\end{proof}

We are now ready to remove the poles as discussed above.
\begin{theorem}\label{thm:coprimefact}
Let $R(\la) \in \F(\la)^{m \times n}$ and $\la_0 \in \F$. There exists a square and regular $K(\la) \in \F[\la]^{n \times n}$ such that
\begin{enumerate}
\item $Y(\la)=R(\la)K(\la) \in \F[\la]^{m \times n}$ does not have poles at $\la_0$;
\item the positive partial multiplicities at $\la_0$ of $Y(\la)$ and $R(\la)$ coincide;
\item the factorization $R(\la)=Y(\la)K(\la)^{-1}$ is right coprime, i.e., the Smith normal form of $\begin{bmatrix}
-K(\la)\\ Y(\la)
\end{bmatrix}$ is $\begin{bmatrix}
I\\ 0
\end{bmatrix}$;
\item if $r(\la)$ is a root vector at $\la_0$ of order $k$ for $Y(\la)$ then $K(\la)r(\la)$ is a root vector at $\la_0$ of order $k$ for $R(\la)$;
\item if $u(\la)$ is a root vector at $\la_0$ of order $k$ for $R(\la)$ then $K(\la)^{-1}u(\la)$ is a root vector at $\la_0$ of order $k$ for $Y(\la)$;
\item the vectors $\{K(\la) r_i(\la)\}$ are a maximal set of root vectors at $\la_0$ for $R(\la)$ if and only if the vectors $\{r_i(\la)\}$ are a maximal set of root vectors at $\la_0$ for $Y(\la)$, with the same orders.
\end{enumerate}
\end{theorem}
It is easy to prove Theorem \ref{thm:coprimefact} using the Smith form. However, the computation of the Smith form of $R(\la)$ (including the unimodular matrices that lead to it) induces the computation of a maximal set of root vectors (they are a certain subset of the columns of one of the unimodular matrices). Instead, we are interested in using Theorem \ref{thm:coprimefact} as a preprocessing tool to compute maximal sets of root vectors when $\la_0$ is a pole (that is, we wish to use it to reduce to the case where $\la_0$ is not a pole). For this reason, we give below a proof which does not rely on the Smith form of $R(\la)$. While this is less direct, it has the advantage of not ``solving the problem in order to solve the problem" in the sense discussed above. Indeed, this proof will inspire our algorithmic approach.
\begin{proof}[Proof of Theorem \ref{thm:coprimefact}]
Without loss of generality, we can represent $R(\la)=D(\la)+C(\la)A(\la)^{-1}B(\la) \in \F(\la)^{m \times n}$ where $A(\la),B(\la),C(\la),D(\la)$ are matrices over $\F[\la]$ of coherent size such that $A(\la)$ is invertible and $P(\la)$ as in \eqref{eq:psm}
is minimal, i.e., the matrices
$ \begin{bmatrix}
-A(\la) & B(\la)
\end{bmatrix}$ and $\begin{bmatrix}
-A(\la)\\
C(\la)
\end{bmatrix}$ have both trivial Smith form over $\F[\la]$.
Then, by \cite[Theorem 3.3]{amz}, there exists a unimodular polynomial matrix $U(\la)$ such that
$\begin{bmatrix}
-A(\la) &
B(\la)
\end{bmatrix} U(\la) = \begin{bmatrix}
I & 0
\end{bmatrix}.$
Hence, there exist polynomial matrices $X(\la)$ and $Y(\la)$ such that
$$P(\la) U(\la) = \begin{bmatrix}
I & 0\\
X(\la) & Y(\la)
\end{bmatrix}$$
implying by Rosenbrock's theorem \cite{C72} that the polynomial matrix $Y(\la)$ has the same positive partial multiplicities at $\la_0$ as $R(\la)$, and (by the fact that it is polynomial) no negative partial multiplicities at $\la_0$.
Similarly, for some rational matrices $W(\la),Z(\la),$ it holds
$$ \begin{bmatrix}-A(\la) & B(\la)\\
0 & R(\la)
\end{bmatrix} U(\la) = \begin{bmatrix}
I & 0\\
W(\la) & Z(\la)
\end{bmatrix}.$$
On the other hand, over $\F(\la)$,
$$ \begin{bmatrix}-A(\la) & B(\la)\\
0 & R(\la)
\end{bmatrix}=  \begin{bmatrix}
I & 0\\
C(\la) A(\la)^{-1} & I
\end{bmatrix} P(\la).$$
Hence, $Z(\la)=Y(\la)$, implying that
$ Y(\la) = \begin{bmatrix}
0 & R(\la)
\end{bmatrix} U(\la) \begin{bmatrix}
0\\
I
\end{bmatrix}.$
In other words, $Y(\la)=R(\la) K(\la)$ for some square  polynomial matrix $K(\la):=\begin{bmatrix}
0 & I
\end{bmatrix}U(\la)\begin{bmatrix}
0\\
I
\end{bmatrix}$. By construction, neither $Y(\la)$ nor $K(\la)$ have poles at $\la_0$, and the positive partial multiplicities at $\la_0$ of $R(\la)$ and $Y(\la)$ coincide. This proves items 1--2. 

Moreover, arguing as in \cite[Theorem 9]{C72}, it is not difficult to prove that $K(\la)$ is regular.
Finally, a slight modification of the arguments in \cite{C72} also shows that
$$Q(\la):=\begin{bmatrix}
-K(\la) & I\\
Y(\la) & 0 
\end{bmatrix}$$
is a minimal system matrix of $R(\la)$ (proving item 3). 

Now suppose that $r(\la)$ is a root vector at $\la_0$ of order $k$ for $Y(\la)$. We claim that $q(\la)=\begin{bmatrix}r(\la)\\
K(\la)r(\la)\end{bmatrix}$ is a root vector at $\la_0$ of order $k$ for $Q(\la)$. Indeed, clearly $\kappa[q(\la)]=\kappa[r(\la)]=0$ and $\kappa[Q(\la)q(\la)]=k$. Moreover, if $q(\la_0) \in \ker_{\la_0} Q(\la)$, then by \cite[Lemma 2.9]{DopN21} there is a polynomial vector $p(\la)$ such that $p(\la_0)=q(\la_0)$ and $Q(\la)p(\la)=0$. Partitioning $p(\la)=\begin{bmatrix}p_1(\la)\\
p_2(\la)\end{bmatrix}$ this in turn implies $p_1(\la_0)=r(\la_0)$ and $Y(\la)p_1(\la)=0$ which, again by \cite[Lemma 2.9]{DopN21}, contradicts the assumption $r(\la_0)\not\in \ker_{\la_0}Y(\la)$. At this point it suffices to observe that $q(\la)$ also satisfies the extra condition in Lemma \ref{lem:sufficient} to prove item 4.

Conversely, assume now that $u(\la)$ is a root vector at $\la_0$ of order $k$
for $R(\la)$. Then, by Theorem \ref{thm:extraction}, $\begin{bmatrix}K(\la)^{-1}u(\la)\\
u(\la) \end{bmatrix}$ is a root vector at $\la_0$ of order $k$ for $Q(\la)$. Defining $v(\la):=K(\la)^{-1} u(\la)$, so that $u(\la)=K(\la)v(\la)$, note that the previous observation implies $\kappa[v(\la)]\geq 0$. Suppose for a contradiction that $v(\la_0) \in \ker_{\la_0} Y(\la)$. Using \cite[Lemma 2.9]{DopN21}, this implies the existence of a polynomial vector $z(\la)$ such that $z(\la_0)=v(\la_0)$ and $Y(\la)z(\la)=0$. Hence, $R(\la)[K(\la) z(\la)]=0$ and $K(\la_0)z(\la_0)=u(\la_0)$, yielding $u(\la_0) \in \ker_{\la_0} R(\la)$: a contradiction. In particular, it cannot be $\kappa[v(\la)]>0$ lest $v(\la_0)=0\in \ker_{\la_0} Y(\la)$, and this concludes the proof of item 5.

For item 6, let us first start from a maximal set $\{r_i(\la)\}_{i=1}^s$ for $Y(\la)$. Similarly to what we did before, we can construct from them a maximal set $\{q_i(\la)\}_{i=1}^s$ for $Q(\la)$ that, moreover, (a) satisfy also the special property of Lemma \ref{lem:sufficient} (b) have bottom block $K(\la)r_i(\la)$. Applying Lemma \ref{lem:sufficient}, we conclude that $\{K(\la)r_i(\la)\}_{i=1}^s$ is a maximal set for $R(\la)$. For the converse implication, let $M(\la)$ be a minimal basis for $\ker R(\la)$. From the analysis in \cite{admz}, the fact that $R(\la)$ and $Y(\la)$ must have the same normal rank, and the fact that $Q(\la)$ is minimal, we see that $N(\la)=K(\la)^{-1}M(\la)$ is a polynomial basis for $\ker Y(\la)$; and if $N(\la_0)$ did not have full column rank, then $M(\la_0)=K(\la_0)N(\la_0)$ also would not: a contradiction. Hence, given a maximal set for $R(\la)$ of the form $\{K(\la)r_i(\la)\}_{i=1}^s$ (any root vector of $R(\la)$ can be written in such form by the previous items),
\[ \rank \begin{bmatrix}M(\la_0) & K(\la_0)r_1(\la_0) & \dots & K(\la_0)r_s(\la_0) \end{bmatrix} \leq \rank \begin{bmatrix}N(\la_0) & r_1(\la_0) & \dots & r_s(\la_0) \end{bmatrix}; \]
since the matrix on the left hand side has full column rank by assumption, so must have the matrix on the right hand side. Hence, $\{r_i(\la)\}_{i=1}^s$ are $\la_0$-independent. As usual, we now invoke Theorem \ref{span} and Theorem \ref{span3}, together with the previously proved fact that $R(\la)$ and $Y(\la)$ have the same positive partial multiplicities at $\la_0$, to also establish completeness and maximality.
\end{proof}

\begin{example}\label{ex:concretK}
The construction of Theorem \ref{thm:coprimefact}, applied to the pair $R(\la),P(\la)$ in Example \ref{ex:bad}, leads to picking (for instance) \[ U(\la)=\begin{bmatrix}0&0&1\\
1&0&\la\\
0&1&0 \end{bmatrix} \Rightarrow K(\la) = \begin{bmatrix}0&\la\\1&0
\end{bmatrix} \Rightarrow Y(\la) = \begin{bmatrix}0&\la\\
1&1
\end{bmatrix}. \]
The vector $\begin{bmatrix}-1&1\end{bmatrix}^T$ is a root vector at $0$ of order $1$ for $Y(\la)$, and we can recover a root vector at $0$ of order $1$ for $R(\la)$ as
\[ K(\la) \begin{bmatrix}
-1\\
1\end{bmatrix} = \begin{bmatrix}
\lambda\\
-1 \end{bmatrix}. \]
\end{example}
\begin{remark}\label{rem:bound}
It is worth noting that left multiplication times $K(\la)$ increased the degree of the (order $1$) root vector from $0$ to $1$. This is remarkable: while, given $\la_0 \in \F$, a root polynomial at $\la_0$ of order $k$ for a \emph{polynomial} matrix can always be assumed to have degree at most $k-1$ \cite[Section 5]{DopN21}, it may be necessary (when the eigenvalue is also a pole) that a root polynomial at $\la_0$ of order $k$ for a \emph{rational} matrix must have degree higher than $k-1$. This is, indeed, shown by Example \ref{ex:concretK}, because it is easy to prove that no polynomial vector of degree $0$ can be a root vector at $0$ for $R(\la)$.

However, let $-m$ be the smallest partial multiplicity at $\la_0$ of $R(\la)$; note that $m>0$ if we assume that $\la_0$ is also a pole. Then, $(\la-\la_0)^m R(\la)$ is polynomial, and it is not difficult to show that $r(\la)$ is a root polynomial of order $k$ at $\la_0$ for $R(\la)$ if and only if it is a root polynomial of order $k+m$ at $\la_0$ for $(\la-\la_0)^m R(\la)$. Hence, it follows that we can assume without loss of generality that a root polynomial of order $k$ at $\la_0$ for a rational matrix $R(\la)$, whose smallest partial multiplicity at $\la_0$ is $-m$, has degree at most $k+m-1$. For instance, in Example \ref{ex:concretK}, we had $k=m=1$ and indeed we constructed a root polynomial of degree $k+m-1=1$. \end{remark}

\begin{remark}\label{rem:coprimefact}
The ingredients needed for the proof of Theorem \ref{thm:coprimefact} can be extended to cover the case where the base ring is any principal ideal domain whose field of fractions is $\F(\la)$. In particular, we can replace $\F[\la]$ with the ring \eqref{eq:localring}, which by construction is a local domain and thus a p.i.d., and reach analogous conclusions following essentially the same arguments.

In particular, in this case, the starting point of the construction is to express $R(\la)=D(\la)+C(\la)A(\la)^{-1}B(\la)$. Then, we build a matrix $P(\la)=\begin{bmatrix}
-A(\la)&B(\la)\\
C(\la)&D(\la)
\end{bmatrix}$ with entries in \eqref{eq:localring}, with $A(\la)$ regular, and minimal at $\la_0$. Minimality at $\la_0$ is defined by asking that both $\begin{bmatrix}
-A(\la) & B(\la)
\end{bmatrix}$ and $\begin{bmatrix}
-A(\la)\\
C(\la)
\end{bmatrix}$ have trivial Smith form over the ring \eqref{eq:localring}. (The Smith form over \eqref{eq:localring} is often referred to as the local Smith form at $\la_0$.) Note that this is a weaker assumption than $P(\la)$ being minimal and polynomial. Moreover, $U(\la)$ is unimodular over \eqref{eq:localring} (i.e. a rational matrix that satisfies the same assumptions of the matrices $A(\la),B(\la)$ in Theorem \ref{thm:PtoAPB}, or in other words is a rational matrix  without poles or zeros at $\la_0$), and $K(\la)$ and $Y(\la)$ are both rational matrices without poles at $\la_0$.

The advantage of working locally is that, instead of removing \emph{all} the poles of $R(\la)$, one would remove only the poles \emph{at the point of interest}. This can be more efficient if one is only interested in the root vectors at one point.

An even further generalization can be realized working over the localization 
\begin{equation}\label{eq:localring2}
{\cal L}
 = \left\{ \frac{a(\la)}{b(\la)} \in \F(\la) : \gcd(a(\la),b(\la))=1 \ \mathrm{and} \ b(\la) \in S \right\}
\end{equation}
 where $S$ is any multiplicatively closed subset of $\F[\la]$ such that $1 \in S, 0 \not\in S$. Indeed, it is not difficult to show that ${\cal L}$ is again a principal ideal domain whose field of fractions is $\F(\la)$.\footnote{Proof: It is a known fact in ring theory that any ideal of ${\cal L}$ has the form $S^{-1}J$ where $J$ is an ideal of $\F[\la]$. But $\F[\la]$ is a p.i.d., and hence $J=\langle g(\la) \rangle$ for some polynomial $g(\la)$. It follows that $S^{-1}J=\langle \frac{g(\la)}{1} \rangle$ is also a principal ideal.} 
\end{remark}

We formalize the analysis of Remark \ref{rem:coprimefact} in the following corollary.

\begin{corollary}\label{cor:coprimefact}
Let $R(\la) \in \F(\la)^{m \times n}$ and $\la_0 \in \F$. Let $S$ be a multiplicatively closed subset of $\F[\la]$ such that $1\in S, 0 \not \in S$, and $(\la-\la_0)^k \not \in S$ for all $k\geq 1$. Moreover, let ${\cal L}$ be the associated ring \eqref{eq:localring2}. There exists a square and regular $n\times n$ matrix $K(\la)$ over \eqref{eq:localring2} such that
\begin{enumerate}
\item $Y(\la)=R(\la) K(\la)$ has elements in \eqref{eq:localring2} and  does not have poles at $\la_0$;
\item the positive partial multiplicities at $\la_0$ of $Y(\la)$ and $R(\la)$ coincide;
\item the factorization $R(\la)=Y(\la)K(\la)^{-1}$ is right coprime (over \eqref{eq:localring2}), i.e., the Smith normal form (over \eqref{eq:localring2}) of $\begin{bmatrix}
-K(\la)\\ Y(\la)
\end{bmatrix}$ is $\begin{bmatrix}
I\\ 0
\end{bmatrix}$;
\item if $r(\la)$ is a root vector at $\la_0$ of order $k$ for $Y(\la)$ then $K(\la)r(\la)$ is a root vector at $\la_0$ of order $k$ for $R(\la)$;
\item if $u(\la)$ is a root vector at $\la_0$ of order $k$ for $R(\la)$ then $K(\la)^{-1}u(\la)$ is a root vector at $\la_0$ of order $k$ for $Y(\la)$;
\item the vectors $\{K(\la) r_i(\la)\}$ are a maximal set of root vectors at $\la_0$ for $R(\la)$ if and only if the vectors $\{r_i(\la)\}$ are a maximal set of root vectors at $\la_0$ for $Y(\la)$, with the same orders.
\end{enumerate}
\end{corollary}

Corollary \ref{cor:coprimefact} suggests a procedure to deal with the case when $\la_0$ is also a pole. The idea is to modify the state space realization of $R(\la)$ to construct a rational matrix $Y(\la)$ as in Remark \ref{rem:coprimefact}, that has no poles at $\la_0$ but has simply related (equal up to premultiplication by $K(\la)$) maximal sets of root vectors at $\la_0$ as $R(\la)$, with the same orders.

In the next subsection, we will show how to realize this in practice by modifying a state space realization of $R(\la)$ in order to build a state space realization of $Y(\la)$.

\begin{remark}\label{manme}
While for many results in the present subsection we have assumed that the polynomial system matrix is (globally) minimal, this assumption can be weakened to local minimality at $\la_0$. We omit a detailed discussion to keep the length of the paper reasonable.
\end{remark}

\begin{remark}\label{rem:new}
Theorem \ref{thm:coprimefact} and Corollary \ref{cor:coprimefact} can in fact be further generalized to another interesting case, by taking the localization ${\cal L}$ to be the ring of proper rational functions
\begin{equation}\label{eq:localring3}
    {\cal L} = \left\{ \frac{a(\lambda)}{b(\lambda)} \in \F(\la) : \gcd(a(\la),b(\la))=1 \ \mathrm{and}  \deg a(\la) \leq \deg b(\la) \right\}
\end{equation} 
instead of \eqref{eq:localring2}. If $\mathcal{L}$ is the ring \eqref{eq:localring3}, we can take $\lambda_0=\infty$ and construct a rational matrix $Y(\lambda)=K(\lambda)R(\lambda)$ with the following properties: (1) $K(\la)$ is square and regular, and the elements of $Y(\la)$ and $K(\la)$ are all proper rational functions (2) $Y(\la)$ has no poles at $\infty$ and has the same positive partial multiplicities at $\infty$ as $R(\lambda)$ (3) the factorization $R(\la)=Y(\la)K(\la)^{-1}$ is coprime over \eqref{eq:localring3} (4) root vectors (resp. maximal sets of root vectors) at $\infty$ for $R(\lambda)$ are $K(\lambda)$ times root vectors (resp. maximal sets of root vectors) at infinity for $Y(\la)$, and vice versa root vectors (resp. maximal sets of root vectors) at $\infty$ for $Y(\la)$ are $K(\la)^{-1}$ times root vectors (resp. maximal sets of root vectors) at $\infty$ for $R(\la)$. (Here, root vectors at infinity should be understood as in the Smith-McMillan approach described in Subsection \ref{sec:infinity}.) 

This claim could be proved again in the same way as Theorem \ref{thm:coprimefact}. For the sake of completeness, we give below a shorter proof via the Smith form. Write $R(\lambda)=U(\lambda)S_R(\lambda)V(\lambda) \in \mathbb{F}(\lambda)^{m \times n}$ where $S_R(\lambda)$ is the $m\times n$ Smith form of $R(\lambda)$ over the ring \eqref{eq:localring3} and $U(\lambda),V(\lambda)$ are unimodular over \eqref{eq:localring3} (i.e. they are square biproper rational matrices). Then without loss of generality $$S_R(\lambda)=\diag(\lambda^{-\sigma_r},\dots,\lambda^{-\sigma_1},0,\dots,0)$$ where $r=\rank R(\lambda)$ and $\sigma_1 \geq \dots \geq \sigma_r$ are integers, and the notation $\diag (m_1,\dots,m_\ell)$ denotes a (possibly rectangular) diagonal matrix with diagonal elements $m_1,\dots,m_\ell$. Let $p$ be the integer such that $\sigma_p \geq  0 > \sigma_{p+1}$ ($0\leq p \leq r$); in other words $r-p$ is the number of poles at $\infty$ of $R(\lambda)$. Define the $n \times n$ square matrix $$S_K(\lambda)=\mathrm{diag}(\lambda^{\sigma_r},\dots,\lambda^{\sigma_{p+1}},1,\dots,1).$$
    Now set $K(\lambda):=V(\lambda)^{-1}S_K(\lambda)$ and $Y(\lambda):=R(\lambda)K(\lambda)$. By construction, $Y(\lambda)$ and $K(\lambda)$ are both proper and morover $Y(\lambda)=U(\lambda)S_Y(\lambda)$ with
    $$S_Y(\lambda)=\mathrm{diag}(1,\dots,1,\lambda^{-\sigma_p},\dots,\lambda^{-\sigma_1},0,\dots,0).$$
    Hence, all the claims made above (considering root vectors at $\lambda_0=\infty$) follow easily.
\end{remark}

\subsection{Coalescent pole/zero: an algorithmic removal of poles that allows us to recover maximal sets of root vectors}

In this subsection we specialize to $\F\subseteq\C$ and look at the computation of root vectors of a general rational matrix $R(\la)\in \C(\la)^{m\times n}$ in the case of a coalescent pole and zero at a finite $\la_0$. We also assume for simplicity that $R(\la)$ has no poles or zeros at $\infty$, and hence $E$ is non-singular. We first recall a result about a so-called feedback transformation applied to the system matrix of an irreducible generalized state-space realization of a given transfer function of McMillan degree $d$~:
\begin{equation}  \label{GSS}
S_Y(\la) := \left[\begin{array}{cc} A+BF -\la E & B \\ C+DF & D \end{array}\right] =  \left[\begin{array}{cc} A-\la E & B \\ C & D \end{array}\right]\left[\begin{array}{cc} I & 0 \\ F & I \end{array}\right] = S_R(\la)\left[\begin{array}{cc} I & 0 \\ F & I \end{array}\right].
\end{equation}
\begin{lemma}[\cite{Wonham79}]\label{lem:wonham}
Given a pair of complex matrices $(\hat A,\hat B)$ of sizes $(d\times d, d\times m)$ and an arbitrary set of complex numbers
$\lambda_1, \ldots, \lambda_n$ then there exists a matrix $F$ such that $\hat A+\hat BF$ has those eigenvalues if and only if the pair $(\hat A,\hat B)$ is controllable.
\end{lemma}
Applying this result to the pair $(\hat A,\hat B) := (E^{-1}A,E^{-1}B)$ shows that this also holds for generalized state-space systems that are irreducible.
The minimality of $S_R(\la)$ implies that for all finite $\la$ the following matrices have full row and column rank $d$, respectively~:
\begin{equation}  \label{min} \rank \left[\begin{array}{cc} A-\la E & B \end{array}\right] = \rank \left[\begin{array}{cc} A-\la E\\ C \end{array}\right]=d. 
\end{equation}  
The new system matrix $S_Y(\la)$ clearly has the same Smith structure at all finite points, but the poles of the corresponding transfer function $Y(\la)$
have changed. This is what we exploit in the next Lemma.
\begin{lemma}\label{lem:finally}The transfer functions $R(\la)$ and $Y(\la)$, of the respective system matrices $S_R(\la)$ and $S_Y(\la)$
given in \eqref{GSS}, are related by an invertible right transformation $K(\la)$ with system matrix $S_K(\la)$
\begin{equation}  \label{cancel} Y(\la)= R(\la)K(\la), \quad S_K(\la) = \left[\begin{array}{cc} A+BF -\la E & B \\ F & I \end{array}\right].
\end{equation}  
If $S_R(\la)$ is minimal, then there always exists a matrix $F$ that assigns the eigenvalues of the pencil 
$(A+BF-\la E)$ to lie in any set $\Lambda \subset \C$. Moreover, if the assigned set $\Lambda$ is disjoint from the poles and zeros of $R(\la)$, then
\begin{enumerate}
\item the inverse transformation $K(\la)^{-1}$ is given by the system matrix $S_{K^{-1}}(\la)=\left[\begin{array}{cc} A -\la E & B \\ -F & I \end{array}\right]$ and the system matrices $S_K(\la)$ and $S_{K^{-1}}(\la)$ are both minimal
\item The system matrix $S_Y(\la)$ is minimal and the Smith-McMillan zeros of $R(\la)$ and $Y(\la)$ are the same
\item the factorization $R(\la)=Y(\la)K(\la)^{-1}$ is a coprime factorization over the ring of functions 
\begin{equation}\label{eq:lambdaring}
{\cal R}_\Lambda :=\left\{ r(\la) \in \C(\la) : |r(\la_0)|< \infty  \;\; \forall \la_0\notin \Lambda \right\},
\end{equation}
which are rational functions that can only have poles in $\Lambda$;
\item the vectors $\{K(\la) r_i(\la)\}_{i=1}^s$ are a maximal set of root vectors at $\la_0$ for $R(\la)$ if and only if the vectors $\{r_i(\la)\}_{i=1}^s$ are a maximal set of root vectors at $\la_0$ for $Y(\la)$, with the same orders.
\end{enumerate}
\end{lemma}
\begin{proof}
We first prove the relations \eqref{cancel}. The following expression is known to be a system matrix for the product 
of the transfer functions with system matrices $S_K(\la)$ and $S_R(\la)$ (see \cite{Van90})~:
$$ \left[\begin{array}{cc|c} A -\la E & 0 & B \\  0 & I & 0 \\ \hline C & 0 & D \end{array}\right] \left[\begin{array}{cc|c}  I & 0 & 0 \\ 0 & A+BF -\la E & B \\ \hline 0 & F & I \end{array}\right] = \left[\begin{array}{cc|c}  A -\la E  & BF & B \\  0 & A+BF -\la E & B \\ \hline C  & DF & D \end{array}\right],
$$ 
which is system equivalent to (after applying simple system equivalence transformations)
$$  \left[\begin{array}{cc|c}  A -\la E  & 0  & 0\\  0 & A+BF -\la E & B \\ \hline C & C+DF & D \end{array}\right],
\quad \mathrm{and} \quad
\left[\begin{array}{c|c} A+BF -\la E & B \\ \hline C+DF & D \end{array}\right].
$$
To prove item 1, we first verify with the same procedure that the product of the transfer functions with system matrices $S_K(\la)$ and $S_{K^{-1}}(\la)$ result in the identity matrix.  The controllability of $S_{K^{-1}}(\la)$ follows from
that of $S_R(\la)$ and the observability of $S_{K^{-1}}(\la)$ follows from the fact that the feedforward term $BF$ of the pair $(A-\la E,F)$ moves all eigenvalues of $A-\la E$ to new locations. Finally, the minimality of $S_K(\la)$ follows from that 
of $S_{K^{-1}}(\la)$, since they both have the same McMillan degree. \\
The system matrix $S_Y(\la)$ may lose minimality only if there is a pole zero cancellation, but the choice of assigned spectrum precludes this, since the poles are the generalized eigenvalues of $(A+BF-\la E)$ and hence are disjoint from the zeros of $S_Y(\la)$.  The rest of item 2 follows trivially from the constant invertible transformation described in \eqref{GSS}, 
between both (minimal) system matrices.\\
Item 3 follows from the fact that both $K(\la)$ and $Y(\la)$ have their poles at the spectrum of $(A+BF-\la E)$,
and hence are bounded outside $\Lambda$.
Moreover, outside $\Lambda$, the compound matrix $\left[ \begin{array}{c} K(\la) \\ -Y(\la)\end{array}\right]$ is bounded and equal to the Schur complement of the system matrix
$$ \left[ \begin{array}{c|c} A+BF-\la E & B \\ \hline C+DF & D \\-F & -I \end{array}\right]= \left[ \begin{array}{c|c} A-\la E & B  \\ \hline C & D \\ 0 & -I \end{array}\right]
\left[ \begin{array}{cc}I & 0\\ F & I \end{array}\right],
$$
which has full column rank because the original system was minimal. 

Finally, for item 4 we observe that maximal sets for $S_R(\la)$ and $S_Y(\la)$ are simply related by Theorem \ref{thm:PtoAPB}. From a maximal set for $R(\la)$, we can construct a maximal set for $S_R(\la)$ via Theorem \ref{thm:extraction}, whence we construct a maximal set for $S_Y(\la)$ which in turn induces a maximal set for $Y(\la)$ by Proposition \ref{prop:conditionalconverse}. From a maximal set for $Y(\la)$, we can go via Theorem \ref{thm:extraction} to a maximal set for $S_Y(\la)$ that, in addition, also \emph{satisfies the extra condition of Lemma \ref{lem:sufficient}}. The constant transformation that relates $S_Y(\la)$ and $S_R(\la)$ is such that this special property is preserved: hence, we can then apply Lemma \ref{lem:sufficient} to such a maximal set for $S_R(\la)$ and recover a maximal set for $R(\la)$. We omit the (straightforward, but tedious) algebraic details necessary to verify that the result of this sequence of transformations eventually lead to the simple relation between maximal sets of $R(\la)$ and $Y(\la)$ that is described in item 4.

\hfill
\end{proof}

\medskip

The analysis of the present subsection can be summarized in the following algorithm:
\begin{itemize}
    \item[Input] A rational matrix $R(\la)$ and a coalescent pole/zero $\la_0$;
    \item[Step 1] Form a minimal system matrix $S_R(\la)$ for the transfer function $R(\la)$;
    \item[Step 2] Find $F$ as in Lemma \ref{lem:wonham} and construct a new system matrix $S_Y(\la)$ as in Lemma \ref{lem:finally};
    \item[Step 3] Compute, for example via the algorithm proposed in \cite{NofV}, a maximal set of root polynomials for $S_Y(\la)$ and denote them in the block form
    $\left\{ \begin{bmatrix}t_i(\la)\\
    b_i(\la) \end{bmatrix} \right\}_{i=1}^s$;
    \item[Step 4(a)] Obtain a maximal set of root vectors for $R(\la)$ as $\{ K(\la)b_i(\la)\}_{i=1}^s$ where $K(\la)$ is as in Lemma \ref{lem:finally}.
\end{itemize}
Step 4(a) requires some additional comment, as $K(\la)$ is generally rational and only available through its representation $K(\la)=I+F(\la E-A-BF)^{-1}B$. Thus, one may wonder how to compute $K(\la) b_i(\la)$ in practice. One option is to obtain $K(\la)$ symbolically, but this may be computationally expensive and anyway, in this paper, we focus on numerical computations. Moreover, depending on the application, one may wish to compute a maximal set root vectors of $R(\la)$ in the special case where they are polynomial vectors, which is not guaranteed by the form $K(\la)b_i(\la)$. To achieve these goals, we can replace Step 4(a) by the following ``truncated" algorithm.
\begin{itemize}
\item[Step 4(b)]
Let $-m$ the the smallest partial multiplicity in $R(\la)$ as in Remark \ref{rem:bound}; note that $m$ can be computed from the staircase form of $A-\la E$ (whose calculation can anyway be useful to reduce the overall complexity: see the discussion after Example \ref{ex:final} below). Then, assuming for notational simplicity $\la_0=0$ (if not, we can perform a linear change of variable), denote $M:=A+BF$ and observe that $M$ is not singular because, by construction, we removed the zero eigenvalues of the pencil $\la E - M$. Hence,
\[ K(\la) = I + F(\la E - M)^{-1} B  = I  - F \left( \sum_{j=0}^\infty \lambda^j (M^{-1}E)^j \right) M^{-1} B. \] 
Suppose now that $b_i(\la)$ is a root polynomial of order $k_i$ at $\la_0=0$ for $Y(\la)$, then using Remark \ref{rem:bound} we can compute a corresponding root polynomial of order $k$ at $\la_0=0$ for $R(\la)$ observing that
\[ K(\la) b_i(\la) \equiv b_i(\la)  - F \left( \sum_{j=0}^{k_i+m-1} \lambda^j (M^{-1}E)^j \right) M^{-1} B b_i(\la) \mod \la^{k_i+m-1}.\]
This operation can be applied to a maximal set of root polynomials for $Y(\la)$ (having possibly distinct orders $k_i$) and produces a maximal set of root polynomials for $R(\la)$.
\end{itemize}

\begin{example}\label{ex:final}
We illustrate the above procedure for Step 4(b) in the case of the rational matrix of Example \ref{ex:bad}. Let us first set $S_R(\la)$ equal to $P(\la)$ in Example \ref{ex:bad} so that 
\[ A=0, \; E=1, \; D=\begin{bmatrix}1&0\\
0&1 \end{bmatrix}, \; B=\begin{bmatrix}1&0 \end{bmatrix}, \; C=\begin{bmatrix}0\\
1\end{bmatrix}.\] We can then pick 
\[ F=\begin{bmatrix}-1\\-1\end{bmatrix} \Rightarrow S_Y(\la)=\begin{bmatrix} -\la-1&1&0\\
-1&1&0\\
0&0&1
\end{bmatrix} \Rightarrow Y(\la)=\begin{bmatrix} \frac{\la}{\la+1}&0\\0&1 \end{bmatrix}. \]
At $\la_0=0$, $S_Y(\la)$ has the maximal set of root polynomials $\begin{bmatrix}
1\\
1\\
0\end{bmatrix}$, whence we extract the maximal set of root polynomials at $0$ for $Y(\la)$ consisting of the single vector $b(\la)=\begin{bmatrix}
1\\
0\end{bmatrix}$. Its order is $k=1$, while the smallest partial multiplicity at $0$ for $R(\la)$ is $-m=-1$. (Note that $m$ is the largest size of a Jordan block with eigenvalue $0$ in the pencil $\la E - A$.) Hence, we can look for a root polynomial at $0$ for $R(\la)$ of degree at most $k+m-1=1$. In particular, $M=A+BF=-1$ and therefore
\[K(\la) b(\la) \equiv \begin{bmatrix}
1\\
0\end{bmatrix} - \begin{bmatrix}-1\\-1\end{bmatrix} \left( \sum_{j=0}^1 (-\la)^j \right) \begin{bmatrix}-1&0 \end{bmatrix} \begin{bmatrix}
1\\
0\end{bmatrix} = \begin{bmatrix}\la\\
\la-1\end{bmatrix}\mod \la.\]
Thus, Step 4(b) constructs $\begin{bmatrix}\la\\
\la-1\end{bmatrix}$ as a maximal set of root polynomials at $0$ for $R(\la)$.
Alternatively, one may apply Step 4(a) and compute exactly
\[ K(\la) = \begin{bmatrix} \frac{\la}{\la+1}&0\\
-\frac{1}{\la+1}& 1
\end{bmatrix} \Rightarrow K(\la)b(\la) = \begin{bmatrix}\frac{\la}{\la+1}\\
-\frac{1}{\la+1}
\end{bmatrix},\]
which is a different maximal set of root vectors for $R(\la)$. In this toy example, Step 4(a) is not difficult, but generally it can be problematic because an exact, say symbolic, computation of $K(\la)$ may be impractical or demanding. Instead, Step 4(b) can always be implemented in a purely numerical manner.
\end{example}

To conclude the paper, we note that the technique described in this subsection allows us to form a rational matrix $Y(\la)$ that has the same zero structure and corresponding root polynomials as $R(\la)$, but has no coalescent pole/zeros anymore. But if we are only interested in one single point $\la_0$ which happens to be also a pole of $R(\la)$, then there exists a more economical procedure. We can then follow essentially
the same idea, but restricted to that pole/zero only.
As pointed out in Section \ref{Sec:systemmatrix}, any irreducible generalized state-space realization 
of $R(\la)$ can be updated to have the following block triangular form
\begin{equation} \label{local0} S(\la) := \left[ \begin{array}{cc|c}
A_1  - \la E_1 & A_2-\la E_2 &  B_1 \\
0 & A_0 - \la E_0 &  B_0 \\ \hline
C_1 & C_0 & D \end{array} \right],
\end{equation}
where the spectrum of the pencil $A_0-\la E_0$ is the single point $\la_0$ and where the pencil $A_1-\la E_1$ 
has no eigenvalues at $\la_0$ and hence contains all the remaining eigenvalues of $A-\la E$. The pole assignment
matrix $F$ is then applied only to the subsystem with pole $\la_0$~:
$$S_Y(\la) =  
 \left[ \begin{array}{cc|c}
A_1  - \la E_1 & A_2-\la E_2 &  B_1 \\
0 & A_0 - \la E_0 &  B_0 \\ \hline
C_1 & C_0 & D \end{array} \right],\left[ \begin{array}{cc|c}
I & 0 & 0 \\
0 &  I & 0  \\ \hline
0 & F_0 & I \end{array} \right],$$
which again corresponds to a factorization 
$$ Y(\la)=R(\la)K(\la), \quad \mathrm{with} \quad K(\la):=F_0(\la E_0-A_0-B_0F_0)^{-1}B_0+I,
$$
and the system matrices
$$
S_K(\la) = \left[\begin{array}{cc} A_0+B_0F_0 -\la E_0 & B_0 \\ F_0 & I \end{array}\right] \quad 
\mathrm{and} \quad S_{K^{-1}}(\la) = \left[\begin{array}{cc} A_0-\la E_0 & B_0 \\ -F_0 & I \end{array}\right].
$$
The minimality of these system matrices follows again from the same arguments, and the factorization 
$R(\la)=Y(\la)K(\la)^{-1}$ is again coprime for very similar reasons. But since we can only guarantee that the 
point $\la_0$ is not a pole anymore of $Y(\la)$, we can only guarantee that the root polynomials of $R(\la)$ 
and $Y(\la)$ are related at $\la_0$. Such a procedure could be repeated for each coalescent pole/zero of $R(\la)$. Moreover, the ``truncated" Step 4(b) can be applied also in this case.

\section{Conclusions} \label{Sec:conclusion}

In this paper we revisited the notion of root polynomials of a polynomial matrix $P(\la)$ at one of its zeros $\la_0$. This notion was introduced in \cite{GLR82} for regular polynomial matrices  and was thoroughly studied in \cite{DopN21} for possibly singular polynomial matrices. An algorithmic approach for root polynomials of pencils was proposed in \cite{NofV}.  We extended the notion of root polynomials to define root vectors of an arbitrary rational matrix $R(\la)$ using valuation theory. This extension makes the definition for coalescent poles and zero consistent with applications. In a second part of the paper, we derive a computational methods to construct a maximal set of root vectors for a general rational matrix $R(\la)$. The approach is based on computing a (locally) minimal state-space realization, applying the algorithm to compute maximal sets of root polynomials of a pencil \cite{NofV}, and finally recovering a maximal set of root vectors for $R(\la)$. If $\la_0$ is a coalescent pole-zero, a preprocessing step implicitly computes a right coprime factorization $R(\la)=Y(\la)K(\la)^{-1}$ such that the maximal sets of root vectors for $Y(\la)$ and $R(\la)$ are simply related. This can be done directly at the level of state-space realizations.

\section*{Acknowledgements}

VN thanks Froil\'{a}n Dopico, Manme Quintana and Ion Zaballa for useful discussions on Rosenbrock's theorem; in particular, credit for some of the very initial steps within the proof of Theorem \ref{thm:coprimefact} should be given to Zaballa who proposed them (in a more general context) while performing research jointly with Dopico and VN. Remark \ref{manme} is due to Quintana, who in addition shared useful comments on a preliminary version of the manuscript. The authors also thank an anonymous referee for useful suggestions and in particular for asking whether Corollary \ref{cor:coprimefact} is also true over the ring of proper rational functions, which led to the addition of Remark \ref{rem:new}.

\end{document}